%% file: TestonMonotonicity.tex
\documentclass[12pt,a4paper,ngerman]{article}
\usepackage{graphicx}
\usepackage{amssymb}
\usepackage{pstricks-add}
\usepackage{pst-plot}
\usepackage{bbm}
\usepackage{inputenc}
\usepackage{ae}
\usepackage{amsmath}
\usepackage{amsthm}
\usepackage{wasysym}
\usepackage{multirow}
\usepackage{subfigure}
\usepackage[round]{natbib}
\usepackage{framed}
\usepackage{makeidx}
\usepackage{array}
\usepackage{verbatim}
\usepackage{pstricks}
\usepackage{longtable}
\usepackage{color}
\usepackage{colortbl}
\usepackage{hyperref}

\makeindex

\newcommand{\bb}{\mathbb }
\newcommand{\vp}{\varphi }
\newcommand{\td}{\tilde }
\newcommand{\ol}{\overline }
\newcommand{\ca}{\mathcal }

\newcommand{\tn}{\, \textnormal}

\newcommand{\R}{\mathbb{R} }
\newcommand{\N}{\mathbb{N} }
\newcommand{\E}{\mathbb{E} }

\newcommand{\beq}{\begin{equation}\label }
\newcommand{\eeq}{\end{equation} }
\newcommand{\bal}{\begin{align*} }
\newcommand{\ve}{\varepsilon}

\newcommand{\bs}{\boldsymbol}
\newcommand{\eal}{\end{align*} }

\newtheorem{satz}{Theorem}[section]
\newtheorem{lemma}[satz]{Lemma}
\newtheorem{defi}[satz]{Definition}

\newtheorem{bem}[satz]{Remark}

\renewcommand{\arraystretch}{1.1}
\makeindex
\begin{document}
\setlength\parindent{0pt}

\parindent 0cm
\title{Multiscale inference for a multivariate density with applications to X-ray astronomy}

\author{
{\small Konstantin Eckle, Nicolai Bissantz, Holger Dette } \\
{\small Ruhr-Universit\"at Bochum } \\
{\small Fakult\"at f\"ur Mathematik } \\
{\small 44780 Bochum, Germany } \\
\and
{\small Katharina Proksch} \\
{\small Institut f\"{u}r Mathematische Stochastik} \\
{\small Georg-August-Universit\"{a}t G\"{o}ttingen } \\
{\small  37077 G\"{o}ttingen, Germany } \\
\and
{\small Sabrina Einecke } \\
{\small Technische Universität Dortmund } \\
{\small Fakult\"at Physik } \\
{\small 44221 Dortmund, Germany } \\
}

\maketitle
\begin{abstract}
In this  paper we propose methods for  inference of the geometric features of a multivariate density.
Our approach uses multiscale tests for the  monotonicity of the density at  arbitrary
points in arbitrary directions.  In particular, a significance  test for a mode at a specific point is  constructed. Moreover, 
we develop multiscale methods for identifying regions of monotonicity
and a general procedure for detecting
the modes of a multivariate density.   It is is shown that the latter method localizes the  modes with
an effectively optimal rate. The theoretical results are illustrated by means of a simulation study
and a data example.
The new method is applied to and motivated by the determination and verification of the position of high-energy sources from X-ray observations by the Swift satellite which
is important for a multiwavelength analysis of objects such as Active Galactic Nuclei.
\end{abstract}

Keywords and Phrases: multiple tests, modes, multivariate density, X-ray astronomy
\\
AMS Subject Classification: 62G07, 62G10, 62G20

\section{Introduction}
\def\theequation{1.\arabic{equation}}
\setcounter{equation}{0}

This work is concerned with the development  of a statistical toolbox which is useful for data
analysis in many problems of  applied sciences. As a specific example we consider a
problem from X-ray astronomy, namely the determination of the positions of objects of
BL Lacertae type with statistical significance. Those objects form a specific subclass of
blazars and hence active Galactic Nuclei (AGN), where high-energy relativistic jets arise
perpendicular to the accretion disc and (in this case) point in the general direction of
the Earth. They are among the brightest extragalactic sources in X-rays and gamma-rays 
in the sky. Determination and verification of the position of such objects from an
observed distribution of origin positions of X-ray photons from the object is of paramount
importance for a multiwavelength analysis of the object to understand and compare the
appearance of photons of different energies in the object. \\
From a statistical point-of-view,  problems of this type are fundamental and refer to a proper understanding 
of the shape of a density $f$ based on a sample of multivariate  observations. Numerous authors
have worked on the detection of qualitative features, such as modes and regions of monotonicity
of a density, in particular on
tests for the existence and  the  localization of  modes. For example, it was pointed out by
\cite{chantong2004} that
the presence of modes can yield to a less precise forecasting. Similarly, a precise
localization of modes  can be used for  non-parametric clustering [see for example
 \cite{MR600539} for an early reference and \cite{MR2332445}, \cite{chaduo2013} and \cite{MR3285753}
 for more recent work].

As pointed out by   \cite{romano1988} and \cite{grundhall1995},  estimation of 
 modes of a density  is a very complex problem, even more difficult than the estimation of the density itself.
In fact, the problem is closely related to the estimation of the first derivative of the
density. There exists a large amount of literature about statistical inference on modes of a density  in the univariate setting, which
 can be roughly divided  into  four different categories: tests on the number of modes, the localization of   modes,
 significance testing of candidate modes and tests that allow for inference about monotonicity.
   \cite{MR1105839} provide the minimax rate  for estimating  a single mode.
 The problem of estimating the  number
 of modes is  considered in  \cite{MR610384} and \cite{hallyork2001}. These authors investigate a test that
 uses bootstrap  methods based on  the so-called critical
bandwidth of a kernel density estimator [see also \cite{mammarfis1992} and \cite{chantong2004} for an asymptotic analysis and
an extension to the dependent case].
In \cite{MR773153},  the distance of the empirical distribution function
to the best-fitting unimodal density is used as test statistic. {\cite{hartigan1987} and  \cite{MR1147099}
propose the excess mass approach for statistical inference of (multi-) modality, which is also used by
\cite{MR1345204,MR1463568} and \cite{MR1844848} to construct  nonparametric tests for the existence of modes.}
\cite{MR1723347} introduce
the SiZer-map as  a graphical tool for the analysis of the local monotonicity properties of a density. 
In this paper, the derivative of a  kernel density estimator is tested locally for a significant deviation  from zero.
 A particular characteristic of the SiZer map is that these tests are performed simultaneously over a fixed range of bandwidths.
A multicsale test for the  monotonicity of a univariate density,  which allows
simultaneous confidence statements about regions of increase and decrease,  can be found in  \cite{MR2435455}. In the univariate
deconvolution model, \cite{MR3113812} propose a multiscale test for qualitative features of a density such as regions 
of monotonicity.
\\
On the other hand,  for  multivariate densities there are just a few results on modality and even less on monotonicity.
\cite{MR1051586} proves
that the optimal minimax rate for mode detection  over a $\beta$-H\"older class is $n^{-{(\beta - 1)}/{(2\beta + d)}}$ $ (\beta >0)$.
{The excess mass approach can also be  used in the multivariate case, but most authors concentrate on one dimension
because - as  pointed out by \cite{MR2508381} -  there is usually a trade-off between practical
feasibility and theoretical justification.}
{\cite{MR2085601}}   use  kernel smoothing  to construct   consistent estimators of the single mode of a multivariate
density, while {\cite{MR2112688}} suggests an adaptive estimate which achieves the optimal rate.
\cite{MR2508381} do not pre-specify the total number of modes
and propose a method for locating   modal regions by means of formal
testing for the presence of anti-modes.  A rate-optimal
algorithm for the localization of the modes of a multivariate density based on a $k$-nearest neighbour
 estimator of the density can be found  in  a recent paper of \cite{NIPS2014_5387}.

A  test about local monotonicity properties of  a bivariate density can be found
 in  \cite{MR1937281} generalizing the SiZer-map. In a multivariate setting \cite{MR2432459}
test locally whether the norm of the gradient of the density vanishes
using a kernel density estimate with  a fixed bandwidth.   In a recent paper
\cite{journals/corr/GenovesePVW13} suggest an algorithm for mode estimation of a $d$-dimensional
density. These authors construct non-parametric confidence intervals for the eigenvalues of the Hessian at
modes of a density estimate, which can be used  for the construction of a significance test.
The method is based on a sample splitting, where   the first half of the
data  is used to localize the modes by means of the mean-shift algorithm
and the second half of the data is used for the significance
test. \cite{journals/corr/GenovesePVW13} also  point  out that the  multiscale approach of  \cite{MR2435455} for constructing confidence
intervals for modes is only applicable to one-dimensional densities.

The  goal of the present paper is to fill this gap
by providing  a multiscale method to identify
 regions of  monotonicity  of a multivariate density.  In Section \ref{sec2} we briefly review the approach
 of    \cite{MR2435455}. We also define a concept of monotonicity in the multivariate case and introduce a multiscale 
 test for this property  at a pre-specified point $x_0 \in \mathbb{R}^d$. The main idea is to investigate monotonicity
 properties of the density in ``various'' directions $e \in \mathbb{R}^d$ by projecting observations from a wedge 
 centered at $x_0$ onto the line $\{ x_0 + te \mid t \geq 0\}$. A multiscale test is provided that allows for 
 a simultaneous inference of the monotonicity properties  at a given confidence level $\alpha$. Section  \ref{sec3}
 extends the approach to the situation where no prior information regarding the location and the number of the modes 
 is available. The theoretical results of this paper establish the consistency of this approach and show that 
  modes can be detected with the optimal rate (up to a logarithmic factor). The finite sample properties of the multiscale 
  test are investigated in Section \ref{sec5} and in Section \ref{secdata} we apply our proposed method to the determination and
  verification of the positions of the blazars Markarian 501 and S3 0218+35.  Finally, all technical details and proofs are deferred to Section \ref{sec6}.

\section{Local testing for a mode}\label{sec2}
\def\theequation{2.\arabic{equation}}
\setcounter{equation}{0}

In this Section, we present a test for the presence of a mode of the density $f$ at a pre-specified candidate point $x_0\in\bb{R}^d$ 
based on a sample of independent random variables   $X_1,\hdots,X_n$ with density $f$. We begin with  a brief review of the work  
of \cite{MR2435455}, who investigate regions of monotonicity of a univariate density (that is $d=1$).

\subsection{Multiscale inference  about a univariate density revisited} 

For one-dimensional  independent identically distributed random variables
$X_1,\hdots,X_n$ with density $f$ let
$X_{(1)}\leq\hdots\leq X_{(n)}$ denote the corresponding  order statistics
and consider  the associated local spacings
\[X_{(i;j,k)}=\frac{X_{(i)}-X_{(j)}}{X_{(k)}-X_{(j)}},\quad j\leq i\leq k.\]
 \cite{MR2435455} propose to use the local spacings
\[T_{jk}(\bs{X})=\sum_{i=j+1}^{k-1}\beta(X_{(i;j,k)}) \]
 to construct a test statistic  for (local) monotonicity of the density $f$ on the interval $(X_{(j)},X_{(k)})$,
  where the function $\beta$ is defined by  $\beta(x):= (2x-1) \mathbbm{1}_{(0,1)}(x) $. Note that $T_{jk}(\bs{X})$ has mean zero if
  $f$ is constant on $(X_{(j)},X_{(k)})$.
Let $X$ denote a random variable with density $f$ independent of $X_1, \ldots, X_n$, and denote by
$$
\tilde F (x) = \mathbb{P} \big(X \leq x \mid X \in [X_{(1)}, X_{(n)}]\big)
$$
the conditional distribution function of $X$ given $X \in [X_{(1)}, X_{(n)}]$.  Define
$U_{(i)} = \td{F} (X_{(i)}) $, then  $U_{(2)} ,\ldots,U_{(n-1)}$ correspond in distribution to the
order statistics of a sample of $(n-2)$ independent uniformly  distributed random variables on the interval $[0,1]$ (note that $U_{(1)}=0$ and $U_{(n)}=1$).
It can be shown that the statistic
\beq{n1}T_{jk}(\bs{U})=\sum_{i=j+1}^{k-1}\beta (U_{(i;j,k)}) \tn{ for }1\leq j<k\leq n,k-j>1,\eeq
satisfies
\[T_{jk}(\bs{X})\begin{cases}\geq T_{jk}(\bs{U}),\tn{ if }f \tn{is increasing on }(X_{(j)},X_{(k)}),\\
                 \leq T_{jk}(\bs{U}),\tn{ if }f \tn{is decreasing on }(X_{(j)},X_{(k)}).
                \end{cases}
\]
Define $
\Gamma(\delta):=\sqrt{2\log ( \tfrac{\exp(1)}{\delta}  )}, $
\beq  {TnU}
T_n(\bs{U})=\max_{1\leq j<k\leq n,k-j>1}\left(\sqrt{\frac{3}{k-j-1}}|T_{jk}(\bs{U})|
-\Gamma\left(\frac{k-j}{n-1}\right)\right),
\eeq
 and  denote by  $\kappa_n(\alpha)$
the $(1-\alpha)$-quantile of  the statistic $T_n(\bs{U})$.  
The multiscale test for  monotonicity proposed by \cite{MR2435455} now concludes that
the density $f$ is not increasing on every
interval $(X_{(j)},X_{(k)})$ with
$$T_{jk}(\bs{X})<-c_{jk}(\alpha) := \sqrt{\frac{k-j-1}{3}}\Big(\Gamma\Big(\frac{k-j}{n-1}\Big)+\kappa_n(\alpha)\Big),\quad 1\leq j<k\leq n,~
k-j>1
$$ and that $f$ is not decreasing on every
interval $(X_{(j)},X_{(k)})$ with $T_{jk}(\bs{X})>c_{jk}(\alpha)$.  The overall risk of at least one
false-positive decision within the simultaneous tests on all scales (i.e. for $ 1\leq j<k\leq n,k-j>1$) is at most $\alpha$.

\subsection{Assumptions and geometrical preparations} 

Throughout this paper $\| x \|$ denotes the Euclidean norm of a vector 
$x \in \mathbb{R}^d$. The  function $f: \R^d \rightarrow \R  $ has a mode at  the point  $x_0$, if for every
vector ${ e}\in\bb{R}^d$ with $\|{ e}\|=1$  the  function  $f_e: t\mapsto f(x_0+t{ e}),\;t\geq 0,$
 is strictly decreasing in a neighbourhood of  $t=0$.
The aim of the test for the presence of a mode defined below is  to investigate the  monotonicity of functions of
this type in different  directions ${ e}$. The number of directions is  determined by the sample size
$n$.
As the set  $\{x_0+t{ e} \mid t\geq0\}$ has Lebesgue measure $0$,
we also have to consider  observations in  a neighbourhood of this line  for inference about monotonicity of the function
 $f_e$. For this purpose, we  introduce a {\it signed distance  of the projection} of a
point $x\in\bb{R}^d$ onto the line
$\{ x_0+t{ e} \mid t\in\bb{R}\}$ and
introduce so-called {\it wedges}. For the following discussion we denote by $\{{ e}_1,\hdots,{ e}_{d-1}\}$ an 
arbitrary but fixed orthonormal basis of $(\mathrm{span}\{ e\})^\perp$, $\langle x,y \rangle$ is the standard inner product of the vectors $x,y \in \mathbb{R}^d$ and
``$\overset{d}{=}$'' denotes equality in distribution.
\begin{defi} {\rm  \label{d1}
Let $x_{0}\in\bb{R}^d$ and $e\in\bb{R}^d$ with $\|e\|=1$.
\begin{itemize}
\item[(1)]   The {\it projected signed distance} of a point $x\in\bb{R}^d$ from $x_0$ in direction ${ e}$ on $\{x_{0}+t{ e} \mid t\in\bb{R}\}$ is defined as
\[P_{ e} x:=\langle x-x_{0},{ e}\rangle.\]
\item[(2)]
The {\it wedge}  with vertex $x_0 $, direction ${ e}$, length $l>0$ and angle  $\varphi\in(0,\frac{\pi}{2})$ is defined as
\bal K\equiv K(x_0,{ e},\vp)&:=\Big\{x\in\bb{R}^d  \ \big  | \; 0<P_{ e}x\leq l\text{ and }\langle x-x_{0},e_i\rangle\in[-\tan(\varphi)P_{ e}x,\tan(\varphi)P_{ e}x]\\&\hspace{0.8cm}\text{ for }i=1,\hdots,d-1\Big\}.\end{align*}
\item[(3)] For a wedge $K \subset \bb{R}^d$ let $X_{(1)} , \ldots , X_{(N)} $ be those random variables among
  $X_1,\hdots,X_n$ which are located in $K$, arranged in ascending order
with respect to their signed projected distances from $x_0$, i.e. $X_{(j)}\in K$ for $j=1,\hdots,N$ 
and $P_eX_{(1)}\leq\hdots\leq
P_eX_{(N)}$.
The wedge $K_N$ is defined as
$K_N:=\big\{x\in K:\; 0<P_ex\leq P_eX_{(N)}\big\}$.
\end{itemize}
}
\end{defi}

\bigskip
\begin{center}
\begin{figure}[ht]
\input{keil.tex} 
\caption{\emph{The wedges $K$ and $K_N$ for $d=2$.}}\label{qw12312312}\end{figure}
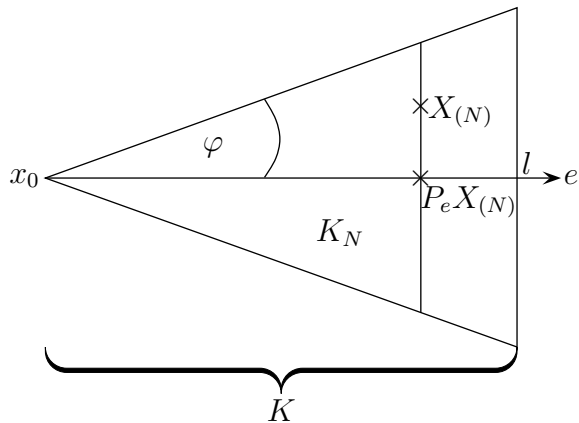
\end{center}\vspace{-1cm}
A typical wedge is  displayed in Figure \ref{qw12312312} in the case $d=2$.
We are now able to   define monotonicity properties of the function $f$ on  the wedge $K$ that will play a crucial role in
following discussion.
\begin{defi} {\rm \label{d2} ~~
\begin{enumerate}
\item The function  $f$ is called {\it increasing on the wedge $K$}, if $ f(x_{0}+\td{t}_2e_0)\geq f(x_{0}+\td{t}_1e_0)$ for all $e_0\in\bb{R}^d$ with $\|e_0\|=1$ and $\td{t}_2>\td{t}_1\geq0$ such that $x_0+\td{t}_2e_0\in K$.
\item The function  $f$ is called   {\it decreasing on the wedge  $K$}, if $ f(x_{0}+\td{t}_2e_0)\leq f(x_{0}+\td{t}_1e_0)$ for all $e_0\in\bb{R}^d$ with $\|e_0\|=1$ and $\td{t}_2>\td{t}_1\geq0$ such that $x_0+\td{t}_2e_0\in K$.
\end{enumerate}
}
\end{defi}

\subsection{A local test for modality}\label{2.2}

  Let  $x_0\in\bb{R}^d$ denote a  candidate position for a mode.
The construction of a local test for the hypothesis that $f$ has a mode at $x_0$ is  based on  an investigation of the
monotonicity properties of $f$ on pairwise disjoint wedges $K^i$   ($i=1,\hdots,M_n$) with common vertex $x_0$.
We  begin with the case $M_n=1$ and use the notation $K:=K^1$ for the sake of simplicity. Throughout this paper $\mathbbm{1}_A$ denotes
the indicator function of a set $A$.

\begin{satz}
\label{P1}
Let $X$ be a $d$-dimensional random variable with density $f$ independent of $X_1,\hdots,X_n$ and denote by
\[ \tilde{F}(z)=\bb{P}(P_eX\leq z|N,X\in K_N,X_{(N)})~~Ê(0<z\leq P_eX_{(N)})
\]
 the  distribution function of $P_eX$ conditional on $N=\sum^n_{i=1} \mathbbm{1}_K(X_i), \{X\in K_N\}$ and $X_{(N)}$.
Then, conditionally on $N$,
$ \td{F}(P_eX_{(1 )}), \ldots ,  \td{F}(P_eX_{(N-1)})$  are distributed  as the order statistics of $N-1$ independent
uniformly distributed random variables on the interval $[0,1]$.
\end{satz}

The first step in the construction of a test  for a mode at the point $x_0$ is to investigate monotonicity in the sense 
of    Definition
\ref{d2}. For this purpose, we use  a comparison of the projected distances $P_eX_{(1)},\hdots,P_eX_{(N-1)}$  with the
distribution
of projected distances  of random variables $U_1,\ldots,U_{N-1}$ which are  uniformly distributed on the wedge $K_N$. For a
random variable $U$ which is uniformly distributed on $K$ and independent of $X_1,\hdots,X_n,U_1,\hdots,U_{N-1}$, we have 
\[\td{F}_U(z)=\bb{P}(P_eU\leq z|N,U\in K_N,X_{(N)})=\frac{z^d}{(P_eX_{(N)})^d} ~~Ê(0<z\leq P_eX_{(N)}),
\]
and by Theorem \ref{P1}, the random variable $\td{F}_U(P_eU_j)=\big(\frac{P_eU_j}{P_eX_{(N)}}\big)^d,~j=1,\hdots,N-1,$ has
 a uniform distribution on the interval $[0,1]$, conditionally on $N$ and the event $\{U_j\in K_N\}$.
Consequently,  we propose  the test statistic
\beq{q10}T_K=\sum_{j=1}^{N-1}\beta\left(\frac{(P_eX_{(j)})^d}{(P_eX_{(N)})^d}\right)\eeq
for testing  monotonicity properties of the density $f$ on the wedge $K$, where  $\beta(z)= (2z-1) \mathbbm{1}_{(0,1)} (z) $.
If   $f$  is constant on $K$, we have
$
\E [T_K] = 0
$
as $\tfrac{(P_eX_{(1)})^d}{(P_eX_{(N)})^d}, \ldots , \tfrac{(P_eX_{(N-1)})^d}{(P_eX_{(N)})^d} $
have the same distribution as an
order statistic of uniformly distributed random variables on the interval $[0,1]$. On the other hand, if $f$ is increasing on
the wedge  $K$, the observations in $K_N$ tend to have large projected distances
from $x_0$, which results in positive values of the test statistic $T_{K}$. Similarly, if $f$ is decreasing on $K$, it is more likely 
that the test statistic is negative.

\begin{satz}\label{P6}
Let $\tilde {F}$ denote the conditional distribution function defined in Theorem \ref{P1}, $T_K$ be defined in \eqref{q10}
and  $T_{K}^U:=\sum_{j=1}^{N-1} \beta(\tilde{F}(P_eX_{(j)}) )$.
\begin{enumerate}
\item
If $f$ is  increasing on $K$, then  $T_{K}^U
\leq T_{K} $ {(a.s.)}  conditionally on $N$.
\item
If  $f$ is  decreasing on $K$, then $
T_K^U\geq T_{K}$ {(a.s.)}  conditionally on $N$.
\end{enumerate}
\label{01}
\end{satz}
By Theorem \ref{P1}, conditionally on $N$,  the statistic $T_{K}^U$  has the same distribution as  the random variable
$\sum_{j=1}^{N-1} \beta (U_{(j)})$, where $U_1, \ldots , U_{N-1}$  are  independent
uniformly  distributed random variables on the interval $[0,1]$. Therefore, Theorem \ref{01} is the key result to obtain
critical values for a multiscale test.

In the  second step, we combine test statistics of the form $T_K$ for different wedges
to construct a test for  a mode at the point $x_0$.  For this purpose, define
\begin{equation}
l_n:=\Big(\frac{\log(n)}{n}\Big)^\frac{1}{d+4}
\label{ln}
\end{equation}
and construct  a family $\ca{K}_n$  of  $M_n$ pairwise disjoint  wedges $K^1 , \ldots, K^{M_n}$ with common vertex $x_0$, 
length  $C_1\log(n)^\frac{d-1}{d+4}l_n$ and angle
$\vp_n:=\frac{C_2}{2}\log(n)^{-1}$ (for
some constants $C_1,C_2>0$) and  by specifying the central directions
$\{{ e}_n^1,\hdots,{ e}_n^{M_n}\}$ as follows
\begin{itemize}
\item[(1)]
Choose a  direction  $ { e}_n^1 $  with $\|{ e}_n^1\|=1$
\item[(2)] If  $ { e}_n^1, \ldots , { e}_n^{i-1} $ have been specified, then - whenever possible -  choose a vector ${ e}_n^i$
with $\|{ e}_n^i\|=1$ such that for some $\varepsilon > 0$
$$
|\mathrm{angle}({ e}_n^i,{ e}_n^j)|\geq (2+\ve)\arctan\big(\sqrt{d-1}\tan(\vp_n)\big)\tn{ for all } j<i.
$$
\end{itemize}
Note that this procedure  does not define $\ca{K}_n$ in a unique way. However, if ${\cal K}_n$ has been fixed for the  central 
directions $e^1_n,\ldots,e^{M_n}_n$, then the following property holds.  For any normalized vector   ${ e}\in\bb{R}^d$ there
exists a direction  ${ e}_n^i$ such that
\[|\mathrm{angle}({ e},{ e}_n^i)|<(2+\ve) \arctan\big(\sqrt{d-1}\tan(\vp_n)\big)=O\big(\log(n)^{-1}\big).\]
This can be seen easily by deriving a contradiction from the opposite assertion using the expansion $\tan(z)=z(1+o(1))=\arctan(z)$ for $z\rightarrow0$.\\
Now define for each $ K^i \in  \ca{K}_n$
 $$N^i := \sum^n_{j=1} \mathbbm{1}_{K^i}(X_j)   $$
 as the number of observations in the wedge $K^i$ and consider the corresponding statistics  $T_{K^i}$
 and  $T_{K^i}^U$ defined in  \eqref{q10} and Theorem \ref{P6}, respectively. An application of
    Theorem \ref{P1}   on each wedge $K^i$ shows that, conditionally on $N^1,\hdots,N^{M_n}$,   \begin{equation} \label{P7}
   T_{K^i}^U\overset{d}{=}\sum_{j=1}^{N^i-1}\beta(U_{(j)}^i) ~~(i=1,\hdots,M_n),
   \end{equation}
 where $\{U_j^i \, | \,  j=1,\hdots,N^i-1,i=1,\hdots,M_n\}$ are  independent   uniformly distributed  random variables
on the interval $[0,1]$. In particular, the statistics  $ T_{K^1}^U , \ldots , T_{K^{M_n} }^U $ are conditionally
independent.

By means of the representation  \eqref{P7}, the quantile  $\td{\kappa}_n(\alpha)$ defined by the condition
\beq{m1}\bb{P}\Bigl(\max_{i=1,\hdots,M_n}\Bigl(\sqrt{\frac{3}{N^i-1}}|T_{K^i}^U|-\Gamma\Bigl(\frac{N^i}{n-1}\Bigr)\Bigr)\leq
\td{\kappa}_n(\alpha)\Big| N^1,\hdots,N^{M_n}\Bigr)= 1-\alpha\eeq
  can be obtained by numerical simulation, as soon as the numbers of observations $N^1,\ldots,N^{M_n}$ in the wedges $K^1,\ldots,K^{M_n}$ have been specified.  We note that a
 calibration by the term $\Gamma(\tfrac{N^i}{n-1})$
 for various scales (i.e. different values of $N^i$) is necessary to show that the quantile $\td{\kappa}_n(\alpha)$
 is asymptotically bounded [see Section \ref{sec6} for details].

  In a third step, we consider on each of the wedges $K^1, \ldots , K^{M_n}$ two hypotheses, that is
\begin{eqnarray}
\label{incr} H^{incr}_{0,i}: ~f  ~ \mbox{ is increasing on } K^i  ~~~~\mbox{ \emph{versus} } ~~~
H^{incr}_{1,i}:~ f  ~\mbox{ is not increasing on }K^i
\end{eqnarray}
($i=1, \ldots ,M_n$), and
\begin{eqnarray}
 \label{decr} H^{decr}_{0,i}: ~f  ~ \mbox{ is decreasing on } K^i  ~~~~\mbox{ \emph{versus} } ~~~
H^{decr}_{1,i}:~ f  ~\mbox{ is not decreasing on }K^i
\end{eqnarray}
($i=1, \ldots ,M_n$),
where the the notation of an increasing (decreasing) function on the wedge $K^i$ is introduced in Definition \ref{d2}.
The $i$th hypothesis in \eqref{incr} is rejected, whenever
\begin{eqnarray}
\label{incrrej}
T_{K^i}<-\td{c}_{K^i}(\alpha)
\end{eqnarray}
($i=1, \ldots ,M_n$), where the quantile $\td{c}_{K^i}(\alpha)$
is defined by
 $\td{c}_{K^i}(\alpha)=\sqrt{\frac{N^i-1}{3}}\big(\td{\kappa}_n(\alpha)+\Gamma (\frac{N^i}{n-1}) \big) $
($i=1,\hdots,M_n$).
 Similarly, the  $i$th hypothesis in \eqref{decr} is rejected, whenever
\begin{eqnarray}
\label{decrrej}
T_{K^i}>\td{c}_{K^i}(\alpha)
\end{eqnarray}
($i=1, \ldots ,M_n$). The final result of this section specifies the error of at least  one false decision among these
$2M_n$  local level $\alpha$-tests on monotonicity.

\begin{satz}\label{P-8}
Assume that all tests  \eqref{incrrej} and \eqref{decrrej} for the hypotheses \eqref{incr} and \eqref{decr}
are performed ($i=1,\hdots,M_n)$.
The  probability of at least one false rejection of any of the tests is  at most $\alpha$.
\end{satz}

We conclude  this section by showing that the results presented so far
can be used to obtain a consistent multiscale test for the hypothesis that the density
$f$ has a mode at a given point $x_0 \in \R^d$.  The test
 decides for the presence of a mode  at  $x_0$ if every test \eqref{incrrej} for the null hypothesis that $f$ is
increasing on the wedge $K^i,i=1,\hdots,M_n,$  rejects the null. Note that in this case
we use the one-sided  quantiles
 $\td{c}_{K^i}'(\alpha)=\sqrt{\frac{N^i-1}{3}}\big(\td{\kappa}'_n(\alpha)+\Gamma (\frac{N^i}{n-1}) \big) $
in  \eqref{incrrej}, where  $\td{\kappa}_n'(\alpha)$ is defined by the condition
\beq{glei1}\bb{P}\Bigl(\max_{i=1,\hdots,M_n}\Bigl(-\sqrt{\frac{3}{N^i-1}}T_{K^i}^U-\Gamma\Bigl(\frac{N^i}{n-1}\Bigr)\Bigr)\leq
\td{\kappa}_n'(\alpha)\Big| N^1,\hdots,N^{M_n}\Bigr)= 1-\alpha.\eeq

\begin{satz}\label{S1A}
Assume that the density  $f$ is  twice continuously differentiable in a neighbourhood of $x_0$
with $f(x_0)\neq 0$, gradient $\nabla f(x_0)=0 $ and a Hessian  $H_f(x_0)$
satisfying $e_0^\top H_f(x_0)e_0\leq -c<0$ for all $e_0\in\bb{R}^d$ with $\|e_0\|=1$.
Consider  the family of wedges  $\ca{K}_n$ defined in Section \ref{2.2} with
constants $C_1,C_2$ satisfying
\beq {c1}
C_1^{d+4}C_2^{d-1}>\frac{4D^2}{c^2}
\frac{f(x_0)}{d+4},
\eeq
where
\beq {bed2}
D=\frac{\sqrt{2}(2d+2)(d+2)}{{\left(1-\frac{d}{d+2}\right)^{(d+2)/ d}
\Big [
1-\frac{d^2}{2(2d+2)^2} \big (-1+  \big \{ 1+4(\frac{2d+2}{d})^2 \big \} ^{1/2} \big) \Big]^{1/2}
}}.
\eeq
Then, all $M_n$ tests defined in \eqref{incrrej} (using the quantiles $\td{c}_{K^i}'(\alpha)$ instead of $\td{c}_{K^i}(\alpha)$
($i=1,\hdots,M_n$))
reject the null hypothesis with asymptotic probability one as $n\rightarrow\infty$.
\end{satz}

Note that the constant  $D$  in \eqref{bed2} depends  only  on the dimension  $d$. Hence,  the
lower bound  on the constants
$C_1$ and $C_2$  is determined by the shape of the modal region (more precisely  the largest eigenvalue
of the Hessian $H_f(x_0)$ at  $x_0$) as well as by the value of the density at the point $x_0$.

\section{Global inference on monotonicity} \label{sec3}
\def\theequation{3.\arabic{equation}}
\setcounter{equation}{0}
In this section  we extend the local inference on modality at a fixed point to the situation where   no specific candidate
position for the mode can be defined in advance. This is particularly important since there exist several applications  where at most
approximate information about the position of  the modes is available. As in the previous section, let $X_1,\ldots,X_n$ denote independent $d$-dimensional random variables with density $f$. The proposed test for the detection of modes proceeds in several steps.

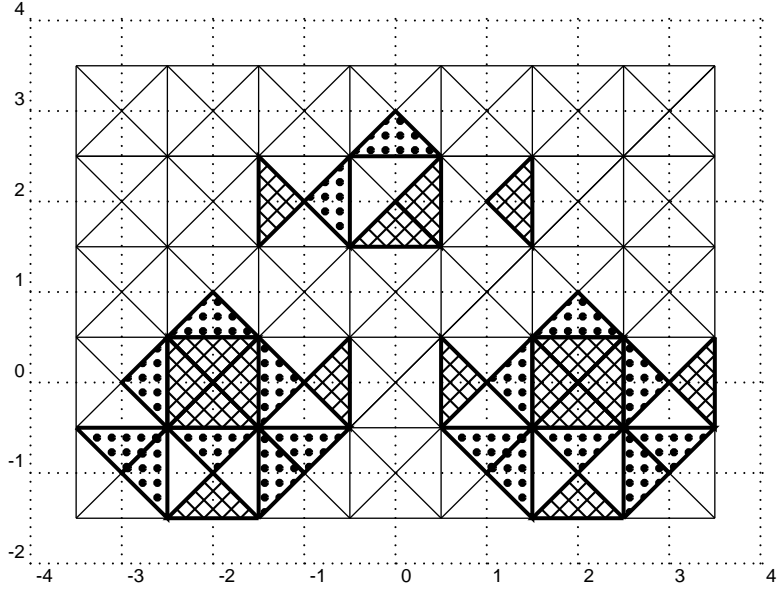
\begin{figure}
\begin{center}
 \input{Gitter_nocolor.tex}
\end{center}
\caption{\emph{Example of a global map for monotonicity of a density.}}\label{f1}
\end{figure}

The first step consists in a  selection of the candidate modes. Here, we choose these as the vertices of an equidistant grid in $\mathbb{R}^d$. Secondly, we introduce
 a  generalization of the multiscale test on monotonicity presented in Section \ref{sec2}, where we divide the wedges in
subsections that are determined by the data. The latter approach can be very useful in settings
without a priori knowledge about the modes, as a true mode obviously has not to be located at the vertex of a  wedge. 
Figure \ref{f1} provides a graphical representation of the results of the global test on modality in the bivariate
case where the multiscale
generalization has been omitted. Here,  on every dotted wedge $K$, the test has rejected that $f$ is decreasing
on $K$. Accordingly,  the cross-hatches refer to a rejection  that $f$ is
increasing on $K$. Non-marked wedges indicate that no significant result has been found.
For a detailed description of the settings used to provide Figure \ref{f1} and an analysis of the  results, we refer
to the end of this section.

\subsection{Geometrical preparations}\label{m3}

Throughout this section let $b_j$ denote  the $j$th unit vector in $\bb{R}^d$  ($j=1,\hdots,d$) and define $\lceil x\rceil:=\inf\{z 
\in\bb{Z}|\,z\geq x\}$.
Recall the definition of $l_n$ in \eqref{ln} and denote (for given constants
 $C_1,\ve>0$)  by $\ca{G}_n$  the grid consisting of the vertices
\[ \sum_{j=1}^d (i_j(2+\ve) C_1\log(n)l_n-n)b_j
\]
 \big($ i_1 \ldots   i_d \in  \big \{0,\ldots,  \lceil\tfrac{2n}{(2+\ve)C_1\log(n)l_n}  \rceil\big \}$ \big).  Note  that the grid $\ca{G}_n$
 covers the cube   $[-n,n]^d$ and   that the sequence $\log(n)l_n$  (which determines the order of the mesh size)
 is chosen such that the test of modality defined below  is consistent.

We now   define by $\ca{K}_n$ a family of wedges (cf. Definition \ref{d1}) with length  $l = C_1\log(n)l_n$, an angle
$\vp = \vp_n=\frac{C_2}{2}\log(n)^{-1}$ for a given constant $C_2>0$,  vertex in $\ca{G}_n$, and a direction
  contained in the set of given directions $\{e_n^1,\hdots,e_n^{M_n}\}$ (cf. Section \ref{2.2}). For
 an arbitrary but fixed element  $K$ of $\ca{K}_n$ let  $X_{(1)}, \ldots , X_{(N)}$ denote  those random variables among
  $\{X_1,\hdots,X_n\}$ which are located  in $K$ and ordered with respect to their signed projected distances from the vertex $x_0$ of $K$.
For  $0\leq j<k\leq N$ with $k-j>1$ we define
\begin{equation*}
K(j,k):=\{x\in\bb{R}^d \ \big  | \;x\in K\tn{ and }P_eX_{(j)}<P_ex\leq P_eX_{(k)}\}
\end{equation*}
as a {\it subsection } of the wedge $K$,
where $X_0:=x_0$ and $x_0$ denotes the vertex of $K$. A typical set is depicted in Figure \ref{fisch}.
We conclude this section with a definition of a
 concept of monotonicity on  subsections of a wedge.

\begin{defi}
Let $K$  be a wedge with vertex $x_0$ and $K(j,k) \subseteq K$ be a subsection. The  function $f: \R^d  \to \R $   is
 \begin{itemize}
  \item[(i)] increasing on $K(j,k)$, if $ f(x_{0}+\td{t}_2e_0)\geq f(x_{0}+\td{t}_1e_0)$ for all $e_0\in\bb{R}^d$ with $\|e_0\|=1$ and
    all $\td{t}_2>\td{t}_1\geq0$ such that $x_0+\td{t}_\ell e_0 \in K(j,k)$ $(\ell=1,2)$.
  \item[(i)]  decreasing on $K(j,k)$, if $ f(x_{0}+\td{t}_2e_0)\leq f(x_{0}+\td{t}_1e_0)$ for all $e_0\in\bb{R}^d$ with $\|e_0\|=1$ and
  all $\td{t}_2>\td{t}_1\geq0$ such that $x_0+\td{t}_\ell e_0 \in K(j,k)$ $(\ell=1,2)$.
 \end{itemize}

\end{defi}

\begin{center}
\begin{figure}[ht]
 \input{Keil2.tex}
\caption{\emph{The subsection   $K(j,k)$ for $d=2$.}}\label{fisch}
\end{figure}
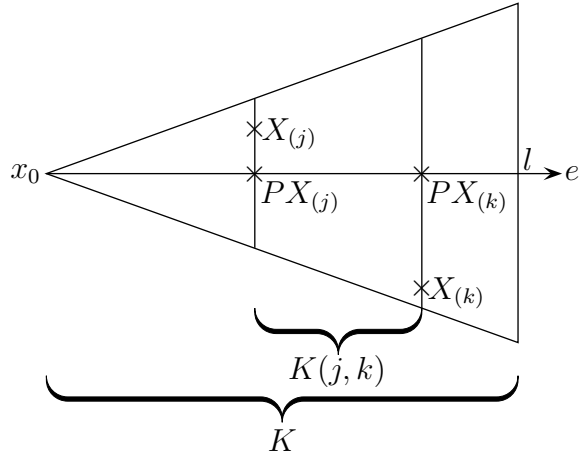
\end{center}

\subsection{Regions of monotonicity and mode detection}\label{3.2}
 The approach proposed here   consists of simultaneous tests for monotonicity of the density $f$ on every subsection
of every wedge in $\ca{K}_n$. For the definition of these tests we will proceed similarly as in Section \ref{2.2}.
We begin by introducing a  multiscale  test statistic on the subsection $K(j,k)$ of a wedge $K\in\ca{K}_n$ which is defined by
\begin{equation} \label{new1}
T_{K(j,k)}:=\sum_{l=j+1}^{k-1}\beta\left(\frac{(P_eX_{(l)})^d-(P_eX_{(j)})^d}{(P_eX_{(k)})^d-(P_eX_{(j)})^d}\right),
\end{equation}
where $0\leq j<k\leq N,\;k-j>1$. Note that $T_{K(0,N)}=T_K$,
 where $T_K$ is the test statistic defined in \eqref{q10}.

Now, let $\ca{K}_n=\{K^i~|\;i=1,\hdots,L_n\}$ denote the family of wedges defined in Section \ref{m3}. 
For the multiscale approach, we use for each subsection
$K^i(j_i,k_i)$   of the wedge $K^i$
the test statistic
$T_{K^i(j_i,k_i)}$ defined by \eqref{new1} ($0\leq j_i<k_i\leq N^i,\;k_i-j_i>1,\;i=1,\hdots, L_n$)
and consider
\[
T_{K^i(j_i,k_i)}^U ~=~
\sum_{l=j_i+1}^{k_i-1} \beta\Big (\frac{\tilde{F}^i(P^i_eX_{(l)})-\tilde{F}^i(P^i_eX_{(j_i)})}{\tilde{F}^i(P^i_eX_{(k_i)})-\tilde{F}^i(P^i_eX_{(j_i)})}\Big )
,\]
 where $P_e^i$ denotes the signed projected distance and
 $\td{F}^i$ denotes the conditional distribution function with respect to 
 $K^i$ ($0\leq j_i<k_i\leq N^i,\;k_i-j_i>1,\;i=1,\hdots, L_n$). Using similar arguments as in Section \ref{2.2},  it follows
that (conditionally on $N^i$)
\begin{eqnarray*}
T_{K^i(j_i,k_i)}^U  & \leq & T_{K^i(j_i,k_i)} ~ (a.s.) ~~\mbox{ if } f \mbox{ is increasing on } K( j_i,k_i),\\
T_{K^i(j_i,k_i)}^U  & \geq & T_{K^i(j_i,k_i)} ~ (a.s.) ~~\mbox{ if } f \mbox{ is decreasing on } K( j_i,k_i)
\end{eqnarray*}
($0\leq j_i<k_i\leq N^i,\;k_i-j_i>1,\;i=1,\ldots, L_n$). Moreover,
\begin{eqnarray*}
T_{K^i(j_i,k_i)}^U &\overset{d}{=}& \sum_{l=j_i+1}^{k_i-1}\beta\Big(\frac{U_{(l)}^i-U_{(j_i)}^i}{U_{(k_i)}^i-U_{(j_i)}^i}
\Big)~~\mbox{conditional on } N^1,\ldots,N^{L_n},
\end{eqnarray*}
 where
 $\{ U_{j_i}^i ~| ~j_i=1,\hdots,N^i-1,i=1,\ldots,L_n \}$ are  independent random variables which are
 uniformly  distributed on the interval $[0,1]$, and
$U_{(1)}^i\leq\ldots\leq U_{(N^i-1)}^i$ is the order statistics of $U_{1}^i,\ldots, U_{N^i-1}^i$
($i=1,\hdots,L_n$).
Finally, let  $\ol{\kappa}_n(\alpha)$ denote the  $(1- \alpha)$-quantile of the conditional distribution
of the random variable
\begin{align}
\label{m2}&
\max_{i=1,\hdots,L_n}\max_{0\leq j_i<k_i\leq N^i,\;k_i-j_i>1}\Big (\sqrt{\frac{3}{k_i-j_i-1}}|T_{K^i(j_i,k_i)}^U|-
\Gamma\Big (\frac{k_i-j_i}{n-1}\Big )\Big )
\end{align}
given $N^1,\hdots,N^{L_n}$. We consider on each subsection $ K^i(j_i,k_i)$ of the wedge $K^i$ the hypotheses
\begin{eqnarray}
\label{incrsub}
H_{0,i,j_i,k_i}^{incr}: f \mbox{ is increasing on } K^i(j_i,k_i) ~  \mbox{versus}   ~ ~H_{1,i,j_i,k_i}^{incr}: f \mbox{ is not increasing on } ~K^i(j_i,k_i),  \\
H_{0,i,j_i,k_i}^{decr}: f \mbox{ is decreasing on } K^i(j_i,k_i) ~\mbox{versus}  ~ ~  H_{1,i,j_i,k_i}^{decr}: f \mbox{ is not decreasing on } ~K^i(j_i,k_i)\nonumber
\end{eqnarray}
($0\leq j_i<k_i\leq N^i,\;k_i-j_i>1$,  $i=1,\hdots,L_n$). The hypothesis $H_{0,i,j_i,k_i}^{incr}$  is rejected if
\begin{eqnarray}
\label{incrtestsub}
T_{K^i(j_i,k_i)}&<& -\ol{c}_{K^i(j_i,k_i)}(\alpha)
\end{eqnarray}
($0\leq j_i<k_i\leq N^i,\;k_i-j_i>1$,  $i=1,\hdots,L_n$), where
$\ol{c}_{K^i(j_i,k_i)}(\alpha):=\sqrt{\frac{k_i-j_i-1}{3}}\big (\ol{\kappa}_n(\alpha)+\Gamma(\frac{k_i-j_i}{n-1})\big).$
Similarly, $H_{0,i,j_i,k_i}^{decr}$ is rejected if
\begin{eqnarray}
\label{decrtestsub}
T_{K^i(j_i,k_i)}>\ol{c}_{K^i(j_i,k_i)}(\alpha)
\end{eqnarray}
($0\leq j_i<k_i\leq N^i,\;k_i-j_i>1$  $i=1,\hdots,L_n$).
Following the line of arguments in the proof of Theorem \ref{P-8}, we obtain the following result.
\begin{satz}
If all tests  \eqref{incrtestsub} and  \eqref{decrtestsub}  are performed
simultaneously $(0\leq j_i<k_i\leq N^i,\;k_i-j_i>1, i=1,\hdots,L_n)$, then
the probability of at least one false rejection  is at most  $\alpha$.
\end{satz}

\subsection{Mode detection} \label{sec33}

We consider the following asymptotic regime. For
$n \in \N$  let ${\cal K}_n$ denote the family of wedges defined in Section \ref{m3} and define
 $\ca{I}_n$  as the set of indices $i$ corresponding to the
wedges  $K^i_n \in {\cal K}_n$ whose vertices  $x_{0,n}$ fulfill $C_1\log(n)l_n\leq\|x_0^n-x_0\|\leq m_nC_1\log(n)l_n
$ for  a mode $x_0$  of $f$ and $m_n=(\log(n))^{\frac{5}{2}}$  and whose direction $e_{n}$ fulfills
$\mathrm{angle}(x_0^n-x_0,e_0^n)=O(\log(n)^{-1}).$
Then, every test \eqref{incrtestsub} for the hypothesis $H_{0,i,0,N^i}^{incr}$
defined \eqref{incrsub} (i.e. $f$ is increasing on $K_n^i$) with
$ i\in \ca{I}_n\subseteq\{1,\hdots,L_n\}$, rejects the null
 with asymptotic probability one.

\begin{satz}\label{S1B}
Let
\[D=\frac{2\sqrt{2}(2d+1)(d+1)}{{\left(1-\frac{d}{d+1}\right)^{(d+1)/d}
\Big[1-\frac{d^2}{2(2d+1)^2}\big(-1+\big\{{1+4(\frac{2d+1}{d})^2\big\}^\frac{1}{2}}\big)}\Big]^{\frac{1}{2}}}.\]

Assume that for any mode  $x_0\in\bb{R}^d$ the density  $f$ satisfies   $c_1\geq f(x_0)>0 $ and
that there exist functions  $g_{x_0}:\bb{R}^d\rightarrow\bb{R}$,   $\td{f}_{x_0}:\bb{R}\rightarrow\bb{R}$
such that  the density $f$
 has a representation of the form
\beq {frep}
 f(x)\equiv (1+g_{x_0}(x))\td{f}_{x_0}(\|x-x_0\|)
\eeq
 (in a neighbourhood of $x_0$).
Furthermore, let $g_{x_0}$ be differentiable in a neighbourhood of $x_0$ with $g_{x_0}(x)=o(1)$ and
$\langle\nabla g_{x_0}(x),e_0\rangle=o(\|x-x_0\|^{1+\gamma})$ (for some $\gamma>0$)
 if   $x\rightarrow x_0$  and all
$e_0\in\bb{R}^d$ with $\|e_0\|=1$. In addition, let $\td{f}_{x_0}$  be differentiable in a neighbourhood
of $0$ with $\td{f}_{x_0}'(h)\leq-ch(1+o(1))$ for  $h\rightarrow0$.
If $\ca{K}_n$ is the family of wedges defined in Section \ref{m3} with
\beq {bed1}
C_1^{d+4}C_2^{d-1} >  \frac {D^2}{c^2}\frac{c_1}{d+4},
\eeq
then
 every mode  $x_0$ of $f$ will be detected with asymptotic probability one as $n\rightarrow\infty.$ 
\end{satz}

Theorem \ref{S1B}  shows that   the proposed procedure  can find all modes with (up to the logarithmic factor) optimal rate.
Note that we proceed in two steps: the verification of the presence of a mode and its localization.
With probability one  the presence of every mode will be detected (by means of the asymptotic
regime introduced at the beginning of this section).
 The rate for the localization of a mode is given by the mesh size of the grid $\ca{G}_n$, which is
determined by the length of the wedges.

\section{Finite sample properties}  \label{sec5}
\def\theequation{4.\arabic{equation}}
\setcounter{equation}{0}
In this section we illustrate the finite sample properties of the proposed multiscale inference. In particular, we 
study the power of the local test for a mode at a given point $x_0 \in \mathbb{R}^d$. We also present
an example illustrating  how the results of Section \ref{sec3} can be used to obtain   a graphical representation  of the 
local monotonicity properties of the density.


\subsection{Local test for modality}

Here, we investigate the  finite sample properties of the local test
for  a two-dimensional 
density, where the level is given by
 $\alpha=0.05$. The corresponding  quantiles $\td{\kappa}_n'(0.05)$ defined in
\eqref{glei1} are determined from $1000$  simulation runs based on
independent and uniformly distributed random variables on the interval $[0,1]$
and are listed in Table \ref{tabquant} for the sample sizes $n=100,500,5000$ in the situation considered in Table \ref{t1} (note that
$\td{\kappa}_n'(0.05)$ depends on the number of observations in every wedge and hence
both on the number and on the size of the wedges).
\begin{table}[h]
 \centering
\begin{tabular}{l| l }
observations&  $\td{\kappa}'_n(0.05)$  \\
\hline
100&0.126\\
500&-0.319\\
5000&-0.854\\
\end{tabular}
\caption{\it Simulated quantiles $\td{\kappa}_n'(0.05)$ in the situation considered in Table \ref{t1}.}\label{tabquant}
\end{table}
 By its construction,
the local test is  conservative, and therefore we also  investigate a calibrated version of the new test.
The quantiles of the calibrated test are chosen such that the level of the test coincides with $\alpha=0.05$ for the
data   obtained from a uniform distribution on the set $[-2.5,2.5]^2$. 
Note that this calibration does not require any knowledge about the unknown density $f$. However, the procedure requires 
the choice of the length and the angle  of the wedges and according to  Theorem \ref{S1A}
we used
\[l_n=C_1 \Big (\frac{\log(n)}{n}\Big)^\frac{1}{d+4}\log(n)^\frac{d-1}{d+4}~,~~\vp_n=\frac{C_2}{2} {\log(n)^{-1}} ,\]
where  $C_1,C_2 >0$ are constants. In the following, the
 power and level of the test   with respect to different choices of  $C_1$ and $C_2$  is investigated.
We also consider
different numbers $M_n$ of wedges in our study.
Recall from the discussion in Section
\ref{2.2} that the constants $C_1$ and $C_2$ have to satisfy \eqref{c1} in order to guarantee consistency of the test.
All results presented below are based on $1000$  simulation runs.

We begin with a comparison of  the test introduced in Section
\ref{sec2} (based  on the  critical values $\td{\kappa}_n^\prime(0.05)$) and a calibrated version of this test.
In Table \ref{t1} we present the simulated level and power of the  local test
for a mode at the point  $x_0=(0,0)^\top$ for different sample sizes.
The constants in the definition of the length and the angle are chosen as  $C_1=2$ and $C_2=9.65$.
For the investigation of the level we consider
a uniform distribution on  the square $[-2.5,2.5]^2$, since it  represents a  ``worst''  case scenario.
For the calculation  of the power, we sample from the standard normal distribution.
We observe that the  test proposed in Section \ref{sec2}  is conservative but it has reasonable power
with increasing sample size.  On the other hand,  the calibrated version of the multiscale test keeps its nominal level  
and rejects the null hypothesis of no mode at $x_0$ in nearly all cases.
\begin{table}[h]
 \centering
\begin{tabular}{l| l | l | l | l | l | l}
observations&  $l_n$ & $M_n$ &level&power&level (cal.)&power (cal.) \\
\hline
100&1.54&3&0.0&36.8&4.8&97.6\\
500&1.31&4&0.0&50.0&4.5&98.4\\
5000&0.99&5&0.0&72.7&5.0&100\\
\end{tabular}
\caption{\it Simulated level and power of the local test for  a mode  at $x_0=(0,0)^\top$ of
 a $2$-dimensional density.}\label{t1}
\end{table}

Next we investigate the influence of the shape of the modal region on the power of the local test. To this end,
we sample from  normal distributions with expectation $(0,0)^\top$ and covariance matrix
$\Sigma\neq I_2$. The results for 
  \renewcommand{\arraystretch}{0.75}
\beq{lp1}\Sigma_1=\Big( \;\begin{matrix} 0&0.5\\-1&1.5 \end{matrix} \;\Big)\tn{ and } ~~ 
\Sigma_2=\Big( \begin{matrix} -0.5&1\\-2&2.5 \end{matrix} \Big) \eeq
are presented in Table \ref{nid}.
\begin{table}[h]
 \centering
\begin{tabular}{l|l|l|l|l|l|l}
\multicolumn{3}{c}{}&\multicolumn{2}{|c|}{$ \Sigma_1$}&\multicolumn{2}{c}{$ \Sigma_2$}\\\cline{4-7}
observations& $l_n$ & $M_n$ &power &power (cal.) &power &power (cal.)  \\
\hline
100&1.54&3&65.4&98.7&38.8&94.3\\
500&1.31&4&95.6&100&80.1&99.6\\
5000&0.99&5&97.8&100&92.1&99.7\\
\end{tabular}
\caption{\it Influence of the shape of the modal region on the power of the local test. The matrices $\Sigma_1$ and $\Sigma_2$ are given in \eqref{lp1}.}\label{nid}
\end{table}
We conclude that the shape of the modal region (determined by the absolute values of the eigenvalues of the
covariance matrix) has a strong influence on the power of the test. In the case $\Sigma_1$ (eigenvalues $0.5$ and $1$), the absolute values
of both eigenvalues are smaller than one.  For $\Sigma_2$ the eigenvalues are given by $0.5$ and $1.5$. Hence,
we observe a (slight) decrease in power in comparison to the first case. However, due to the existence  of an eigenvalue with absolute value smaller
than one,   the test still performs better as in the case of a standard normal distribution.

As the local test requires the specification of the point $x_0$, we next investigate the  influence of
its incorrect specification  on  the power of the test.
For this purpose, we consider the same data (two-dimensional standard normal distribution)
and  perform the tests  under the assumption that the modes are given by  $x_0=(0.2,0.2)^\top$
and  $x_0=(0.7,0.7)^\top$, respectively (which has to be compared to the
true position of the mode at $(0,0)^\top$). The corresponding results are shown in Table \ref{t2},
and we conclude that a ``small''  deviation of the candidate mode from  the true mode has a very
small effect on the power of the tests. In the case $x_0=(0.7,0.7)^\top$, the
 distance between the candidate and the true mode  is very large  in comparison to the length of the wedges. 
 For $n=100$ observations the  length of the wedges is still substantially larger than the distance between the candidate and the true mode. Hence, the test  detects
 the presence of a mode, but we observe a decrease in its
 power. However, for $n=5000 $ observations the distance between the
candidate position and the true mode is approximately equal to the length of the
 wedges. As a consequence, the multiscale test is performed with a finer triangulation and (correctly) 
 does not  indicate the existence
 of  a mode at the point $x_0=(0.7,0.7)^\top$.
\begin{table}[h]
 \centering
\begin{tabular}{l|l|l|l|l|l|l}
\multicolumn{3}{c}{}&\multicolumn{2}{|c|}{$x_0=(0.2,0.2)^\top$}&\multicolumn{2}{c}{$x_0=(0.7,0.7)^\top$}\\\cline{4-7}
observations& $l_n$ & $M_n$ &power &power (cal.) &power &power (cal.)  \\
\hline
100&1.54&3&32.4&96.6&2.8&75.6\\
500&1.31&4&43.1&97.8&1.2&57.1\\
5000&0.99&5&47.2&98.3&0.1&10.8\\
\end{tabular}
\caption{\it Influence of a misspecification of the mode on the power of the local test.}\label{t2}
\end{table}

In the remaining part of this section  we investigate the influence of the choice of the parameters $C_1$ and $C_2$ on the power
of the test. Note that the volume of every wedge is proportional to $l_n^d$, where $l_n$ is
the length of the wedge. This means that dividing the length in half
yields  a wedge with a volume which is $2^{-d}$ times  smaller than the volume of the original wedge.
Thus, the number of observations in the smaller wedge is substantially smaller
 than the number of observations in the larger wedge. Therefore, we expect that the
 constant $C_1$ has  an impact  on the power of the test. These theoretical considerations are reflected by the
 numerical results in
Table \ref{t3}, which show the power for a fixed sample size $n=500$, different choices of $C_1$
(represented by the different lengths) and fixed  parameter $C_2$.  We observe a loss of power of both
tests with decreasing length of the wedge.
 \\
On the other hand,  decreasing the constant $C_2$ such that the number of wedges doubles has the effect that the number
of observations in every wedge decreases approximately by $50\%$.  In Table \ref{t4} we show  the
 power for a fixed sample size $n=500$, a fixed constant  $C_1= 2$ and varying values of $C_2$ (represented by the different number of wedges).
 Here, the picture is not so clear. While we observe a loss in power of the non-calibrated tests with an increasing number of wedges,
 the power of the calibrated test changes only slightly.
 In both cases, the calibrated version still performs rather well  opposite to
its uncalibrated version.

\begin{table}[h]
 \centering
\begin{tabular}{l|l|l|l|l|l|l}
observations& $l_n$ &  $M_n$ &level&power&level (cal.)&power (cal.)  \\
\hline
500&1.31&4&0.0&50.0&4.5&98.4\\
500&0.98&4&0.0&1.5&4.9&74.4\\
500&0.65&4&0.0&0.0&5.4&37.2\\
\end{tabular}
\caption{\it Influence of the length of the wedges on the power.}\label{t3}
\end{table}

\begin{table}[h]
 \centering
\begin{tabular}{l|l|l|l|l|l|l}
observations&  $l_n$ &   $M_n$ &level&power&level (cal.)&power (cal.)  \\
\hline
500&1.31&4&0.0&50.0&4.5&98.4\\
500&1.31&6&0.0&0.7&5.3&92.7\\
500&1.31&8&0.0&0.0&4.9&89.8\\
\end{tabular}
\caption{\it Influence of the number of directions tested on the power.}\label{t4}
\end{table}

\subsection{Identifying local monotonicity of a multivariate density} \label{sec52}
\begin{figure}
\begin{center}
\includegraphics[width=0.6\textwidth]{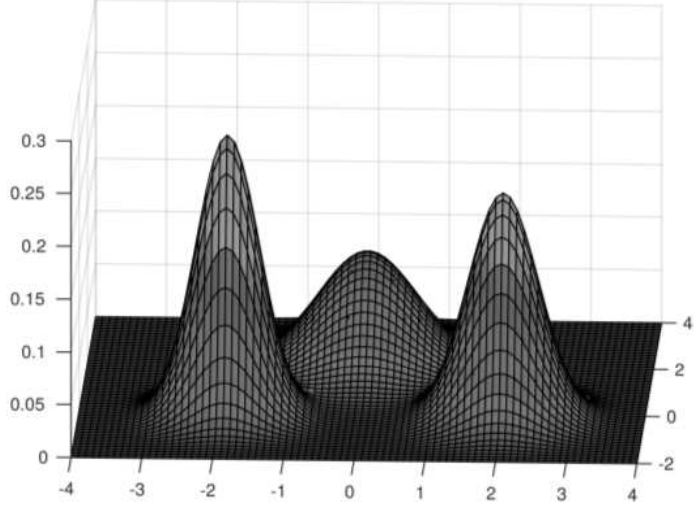}
\end{center}
\caption{\it The density of  a (uniform)   mixture  of  a ${\cal N}((-0.05,2.1)^\top,0.5I),$
${\cal N}((-1.9,-0.07)^\top,0.2I)$ and $ {\cal N}((2,-0.1)^\top,0.25I)$ distribution.} \label{f2}
\end{figure}

In this section we demonstrate how the results of  Section \ref{sec3}
can be used to obtain a graphical representation of the
local monotonicity behaviour of the density (in the case $d=2$).
We conduct the procedure to detect  regions of monotonicity as proposed in Section \ref{3.2}. For the sake of
convenience, we use only
the largest scales in the test statistic \eqref{new1} (i.e. we test on the entire wedges  and not  on the subsections introduced in Section \ref{m3}).
 The significance level  is $\alpha=0.05$.
 We chose an equidistant grid  covering  $[-3,3]\times[-1,3]$ with points $(i,j)^\top,\;i=-3,\hdots,3,j=-1,\hdots,3$, the
 length of any  wedge is   $l=\frac{1}{2}$ and all angles are given by $\vp=\frac{\pi}{4}$.
 Figure \ref{f1} presents the map of the local monotonicity properties  on the basis of $n=100000$
observations from a mixture of three normal distributions  (i.e. $f$  has three modes of different shape)  [see  Figure \ref{f2}].
Here, a cross-hatched wedge indicates that the local  test rejected the hypothesis that  the density is increasing on the respective wedge.
Similarly, a dotted wedge implies that the test rejected the hypothesis that  the density is decreasing. Non-significant
wedges are not marked. 

\begin{figure}[h]
\begin{center}
 \input{Mode2_nocolor.tex}
\end{center}
\caption{\emph{Indications for the presence of a mode. Left panel: grid points $(-2,0)^\top$ and $(2,0)^\top$. 
Right panel: grid point $(0,2)^\top$}.}\label{f4}
\end{figure}
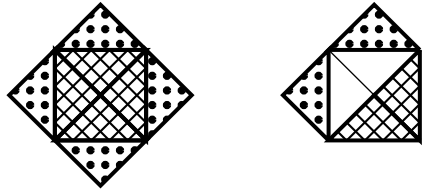

The map  indicates the existence of modes close to the grid points $(-2,0)^\top$ and $(2,0)^\top$ and in a weaker sense
indicates  a mode close to
the grid point at $(0,2)^\top$. The marked geometrical objects around these grid points are shown in Figure \ref{f4}.
In the grid point at $(0,2)^\top$ we obtain not so many significant rejections as in the wedges with vertex
$(2,0)^\top$. Still, the dotted wedges show that there is a significant increase towards the mode which gives an indication
for the presence of a mode as well. An improved procedure with a direct focus on the modes will be discussed in the following section.

\subsection{Mode detection}  \label{sec53}

In this section we demonstrate how the multiscale test can be  successfully used for the localization of modes
if the inference on the local monotonicity behaviour of the density is not included in the  test statistic.
More precisely,  we  consider the grid  introduced in Section \ref{sec52}.
Similar to the local test on modality, we conclude that the density has a mode close to a grid
point if every test on every wedge whose vertex is given by the
grid point rejects that the density is increasing on the respective wedge.
We again recommend  a calibrated version of the global test, where
the  quantiles are chosen such that the probability of the discovery of a
non-existing mode is approximately $5\%$ if the data comes from a two-dimensional uniform distribution.  The following results are based on
$n=2500$ observations and $1000$ simulation runs.  \\
We have investigated two densities, a constant
density on the square  $[-3.5,3.5]\times[-1.5,3.5]$ and  the  density with three modes  presented in  Figure \ref{f2}. For the uniform
distribution  the test  found a mode in $4.6\%$  of the simulation runs. For the
the tri-modal density the test detected  in $78.9\%$  cases  a mode in the point
$(-2,0)^\top$, in $53.8\%$  cases a mode  in $(2,0)^\top$ and in $7.4\%$
cases a mode in the point  $(0,2)^\top$.

\section{Real data example} \label{secdata}
Active Galactic Nuclei (AGN) consist of a supermassive black hole at their centre surrounded by an accretion disc. 
For some AGN classes, high-energy relativistic jets perpendicular to the disc are produced.
In case this jet is pointing in the general direction of the Earth, the object is referred to as a blazar.
The BL Lacertae type object Markarian (Mrk) 501 is such a blazar and one of the closest (in a distance of
$4.7 \times 10^8$ light years) and brightest extragalactic sources in X-rays and gamma-rays.
It is known as one of the most extreme blazars and features very strong and fast variability, making Mrk 501 a perfect 
candidate for probing AGN.
Due to the strong emission over the entire electromagnetic spectrum, correlation studies between different energy bands 
(parts of the electromagnetic spectrum) are particularly interesting and will give insights into the processes inside an
AGN, e.g. the emission models or the particle populations, since different spatial regions of the object may emit 
radiation of diverse energies.
Therefore, the determination of the position of the radiation in a certain energy regime is of paramount importance.

\begin{figure}
\begin{center}
\includegraphics[width=0.5\textwidth]{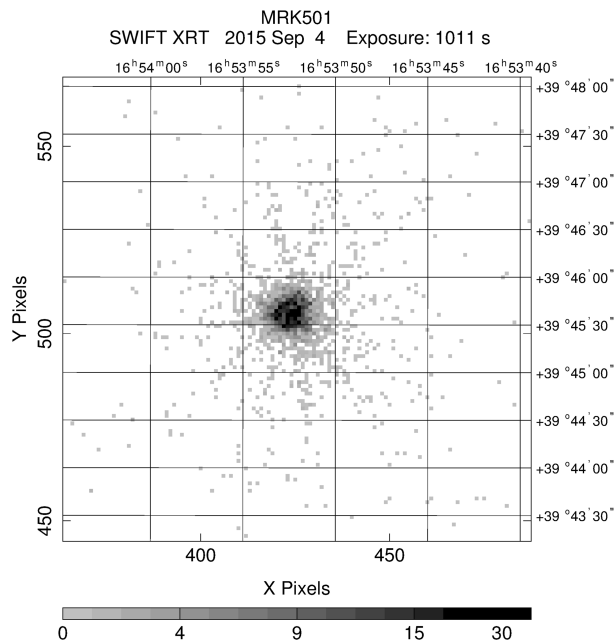}
\end{center}\vspace{-0.8cm}
\caption{\it Observation of Mrk 501 from Sep. 4, 2015.\protect\footnotemark} \label{stern}
\end{figure}

Here, 19 individual observations of Mrk 501 in the year 2015, performed in photon counting (PC) mode by the \textit{Swift}-XRT 
on board the \textit{Swift} satellite, are analysed. 
The \textit{Swift} satellite was launched in 2004 and is a multiwavelength space observatory with three instruments on board, 
the X-ray Telescope (XRT) being one of them, which is capable to observe X-rays in the $0.3-10\,\textnormal{keV}$ energy regime. 
For each considered observation, the High Energy Astrophysics Science Archive Research Center
(HEASARC\footnote{NASA/Goddard Space Flight Center, https://heasarc.gsfc.nasa.gov/.}) provides an image, based on Level 2 event 
files that have been calibrated and screened by a standard pipeline.
These images with a size of $1000 \times 1000$ pixels contain the information how many photons (i.e. which X-ray flux) have been 
recorded in each pixel during the exposure time.
The exposure times of the analysed images range from about $100\,\textnormal{s}$ to $1000\,\textnormal{s}$.
Due to different positions of the satellite in space and different alignments of its main axis, each image shows a
slightly different region of the sky. Figure \ref{stern} provides an illustration of the data obtained from one observation.

\begin{figure}[ht]
\begin{center}
\includegraphics[width=0.4\textwidth]{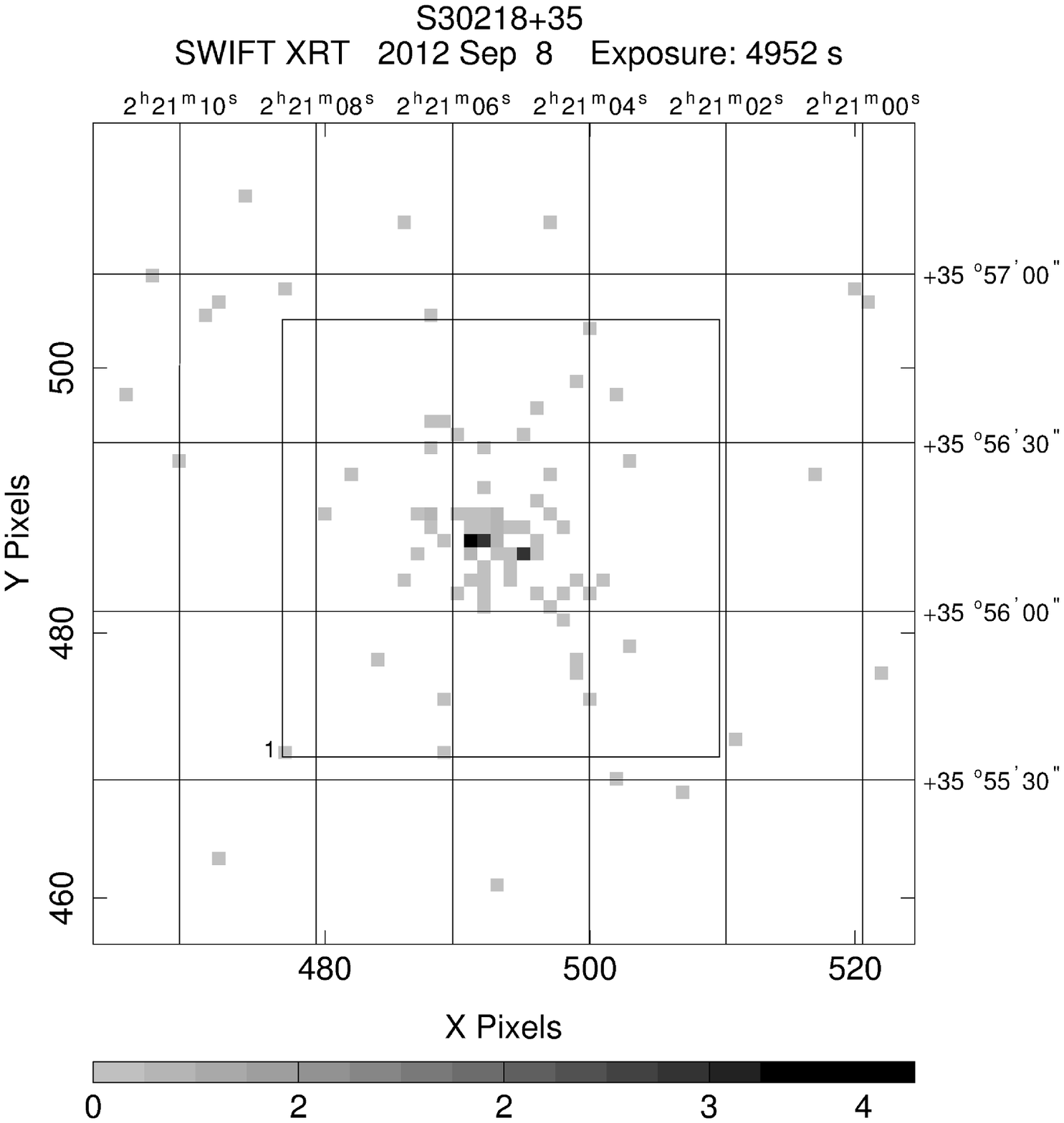}
\includegraphics[width=0.4\textwidth]{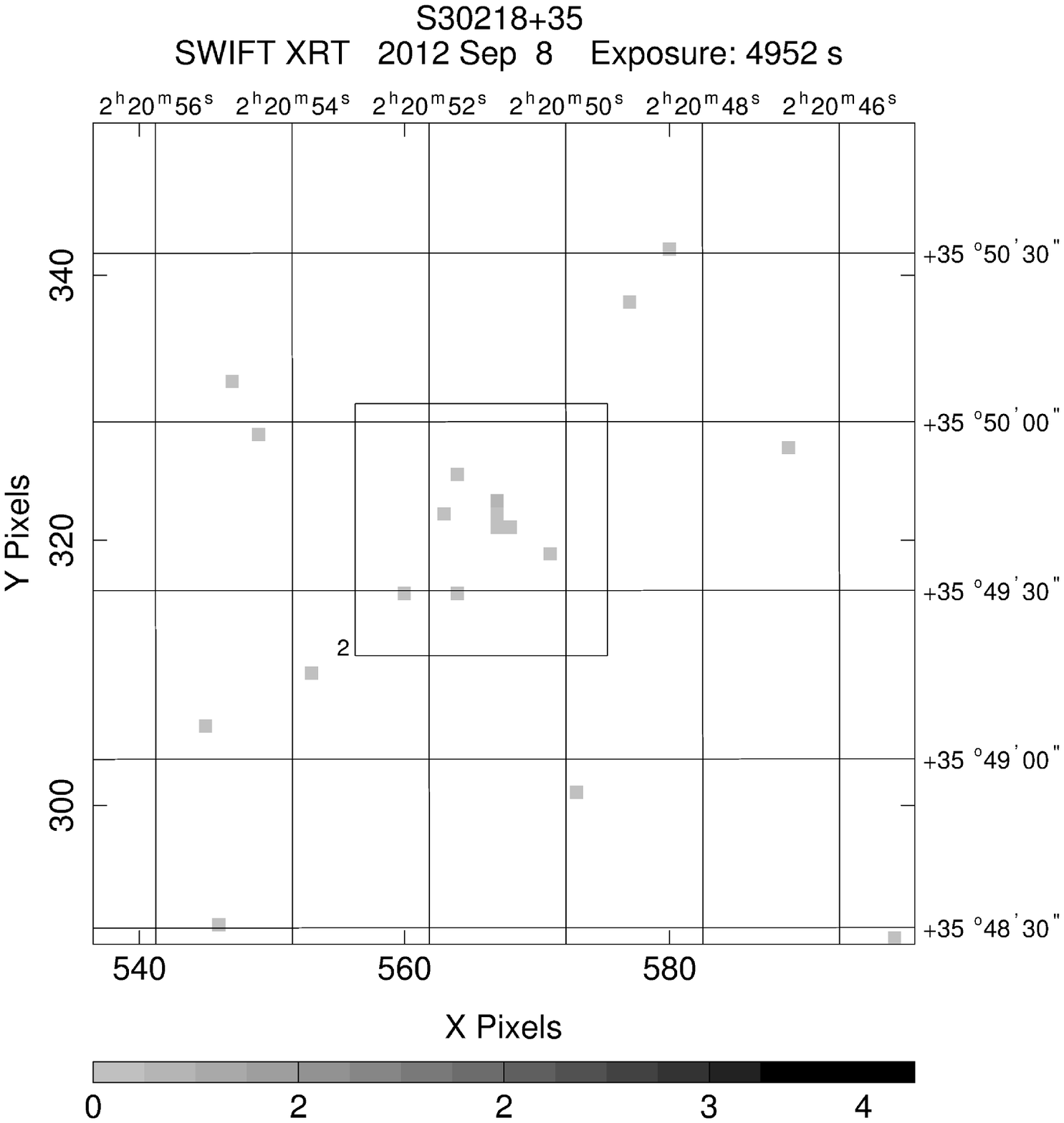}
\end{center}\vspace{-0.8cm}
\caption{\it Observations of the two points sources around S3 0218+35 from Sep. 8, 2012.\protect\footnotemark} \label{z1weistern}
\end{figure}

By a combination of the 19 individual observations, we have at our disposal $49248$ observations of X-ray photons with known positions of
origin on the sky. Our 
aim is the precise localization of the mode of the distribution. To this end, we conduct the test presented in Section
\ref{sec53} for a significance level of $0.01$. We chose an equidistant grid covering 
$[253.446^\circ,253.586^\circ]\times[39.64^\circ,39.88^\circ]$ (corresponding to $[16\tn{h } 53\tn{min } 47.04\tn{s}, 16\tn{h } 54\tn{min } 20.64\tn{s}]
\times  [39^\circ\; 38'\; 24'', 39^\circ\; 52' \;48'']$) consisting of 961
grid points with mesh size $0.008^\circ$. The length of any wedge is  $l=0.004^\circ$ and all angles are given by $\vp=\frac{\pi}{4}$.
Again, we used a calibrated version of the test where the quantiles are chosen such that a non-existing mode for a
uniform distribution on $[253^\circ,253.8^\circ]\times[39.5^\circ,40.1^\circ]$ has been found in less than
$1\%$ of the simulation runs (based on 1000 simulation runs). Our test detected the mode of the distribution in 
$(253.466^\circ,39.760^\circ)$ (corresponding to $(16\tn{h } 53\tn{min	}51.84\tn{s}, 39^\circ	\;45'\;	36'')$). The precision regarding the location of this mode is given by the mesh size $0.008^\circ$.\\
\footnotetext{This picture has been created using HEAsoft, http://heasarc.nasa.gov/docs/software/lheasoft/.}
\begin{figure}[ht]
\begin{center}
\includegraphics[width=0.5\textwidth]{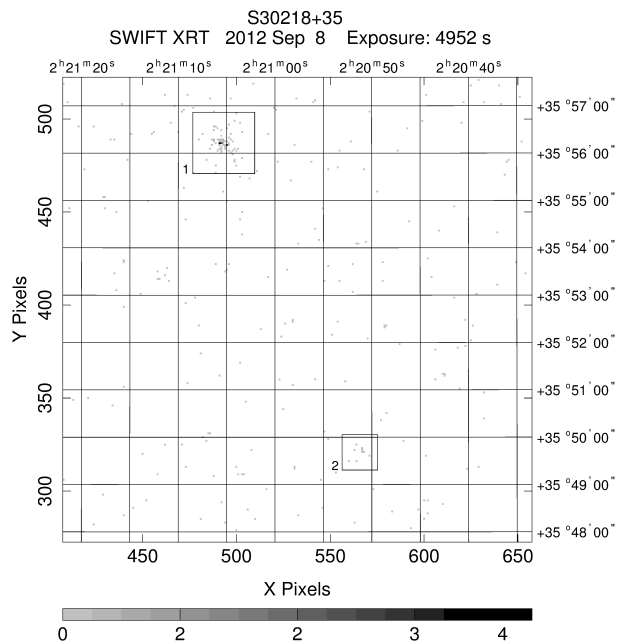}
\end{center}\vspace{-0.8cm}
\caption{\it Observation of the sky region around the blazar S3 0218+35 from Sep. 8, 2012.\protect\footnotemark} \label{zweistern}
\end{figure}
In contrast to the high confidence detection possible within short exposure times for the brightest extragalactic sources, such as Mrk 501, the detection of fainter point sources is more challenging.
Often only few photons reach the detector due to, e.g., the large distance to the source and the absorption of photons.
Within the sky region of one \textit{Swift}-XRT image, there might be multiple point sources in X-rays, but even background 
fluctuations can look like faint point sources.
The study of the population of these point sources, the correlation to other energy bands, and variability studies 
contribute enormously to the understanding of the X-ray sky.
This requires
 reliable methods for the detection and the determination of the position, including the confidence of a given calculation.
In the following, the capability to determine multiple modes of faint point sources in \textit{Swift}-XRT images is demonstrated.
18 images provided by HEASARC of  individual observations of the sky region around the blazar S3 0218+35 in the years 2005, 2012, 2014 and 2016 performed in PC mode by the \textit{Swift}-XRT are analysed.
The exposure times of the images provided by HEASARC range from 3000\,s to 5000\,s.
Figure \ref{zweistern} shows one of these images. The two point sources are marked with a square. Figure \ref{z1weistern} provides
 detailed images of the two point sources.\\\footnotetext{This picture has been created using HEAsoft, http://heasarc.nasa.gov/docs/software/lheasoft/.}

For the detection and the localization of the two point sources we conduct the test presented in Section \ref{sec53} at a significance level of
$0.01$. Here, we have at our disposal 18061 observations. For this application, we chose an equidistant grid covering $[35.2^\circ,35.32^\circ]
\times[35.825^\circ,35.945^\circ]$ (corresponding to $[2\tn{h } 20\tn{min }48\tn{s}, 2\tn{h } 21\tn{min }16.8\tn{s}] \times 
[35^\circ\; 49' \;30'', 35^\circ \;56' \;42'']$)
consisting of 961 grid points with mesh size $0.004^\circ$. The length of any wedge is  $l=0.002^\circ$ and all angles are given by 
$\vp=\frac{\pi}{4}$. The quantiles are chosen such that a non-existing mode for a 
uniform distribution on $[35^\circ,35.5^\circ]\times[35.7^\circ,36.1^\circ]$ has been found in less than
$1\%$ of the simulation runs (based on 1000 simulation runs). Our test detected the two point sources at $(35.212^\circ,35.829^\circ)$
(corresponding to $(2\tn{h }20\tn{min }50.88\tn{s}, 35^\circ \;49'\; 44.4'')$)
and at $(35.272^\circ,35.937^\circ)$ (corresponding to $(2\tn{h }21\tn{min } 5.28\tn{s}, 35^\circ \;56' \;13.2'')$) at a precision of $0.004^\circ$.
 \bigskip
 \bigskip

{\bf Acknowledgements} This research has made use of data obtained through the High Energy Astrophysics Science Archive Research 
Center Online Service, 
provided by the NASA/Goddard Space Flight Center. The authors would like to thank Martina
Stein, who typed parts of this manuscript with considerable
technical expertise.
This work has been supported in part by the
Collaborative Research Center ``Statistical modeling of nonlinear
dynamic processes'' (SFB 823, Teilprojekt C1, C4) of the German Research Foundation
(DFG).
\footnotetext{These pictures have been created using HEAsoft, http://heasarc.nasa.gov/docs/software/lheasoft/.}

\begin{small}
 \setlength{\bibsep}{2pt}
\bibliography{lit3}
\end{small}

\section{Proofs} \label{sec6}
\def\theequation{5.\arabic{equation}}
\setcounter{equation}{0}
The assertions of most of our results are split up in two parts, one for densities that are increasing and one
for decreasing densities. Often, the proof for  one case can be transferred in a straightforward way to the other one. In this 
situation, we only prove one case  as the other case follows immediately by similar arguments.
Only if this transfer is not obvious,  we give   details for both situations.

\subsection{Proofs of  the results of Section \ref{2.2}}
\begin{proof}[\textbf{Proof of Theorem \ref{P1}}] \relax
It is well-known that,  given $X_{(N)}$ and $N$, the random variables $\td{F}(P_eX_{(1)}), \ldots, \td{F}(P_eX_{(N-1)})$
have the same distribution as the order statistics of $N-1$ uniformly  distributed and independent random variables on
the interval $[0,1]$. By an application of the law of iterated expectations,   the assertion follows.
\end{proof}

\begin{proof}[\textbf{Proof of Theorem \ref{01}}]
 We only consider the case where the density $f$ is increasing on $K$. 
 As   $\td{F}(P_eX_{(j)})\neq0$ almost surely for $j=1,\hdots,N-1$ (cf. Theorem \ref{P1}) and $P_eX_{(j)}=P_eX_{(N)}$ 
 implies $\td{F}(P_eX_{(j)})=1$,
it is sufficient to prove
\beq{-15}\td{F}(z)
\leq \td{F}_U(z) =
\frac{z^d}{(P_eX_{(N)})^d}\eeq
for all $z\in(0,P_eX_{(N)}]$.
For this purpose, notice that the distribution function $\td{F}$  in
Theorem  \ref{P1}  is given  by
\bal& \tilde{F}(z)=\frac{\bb{P}(P_eX\leq z,X\in K_N| N,X_{(N)})}{\bb{P}(X\in K_N|N,X_{(N)})} =
\frac{G(z)}{G(P_eX_{(N)})}~,
\end{align*}
where the function $G$ is defined by
\begin{equation*}
G(z) = \int_{0}^z\int_
{-t\tan(\varphi)}^{t\tan(\varphi)}\hdots\int_{-t\tan(\vp)}^{t\tan(\vp)}f(x_{0}+te+s_1e_1+\hdots+s_{d-1}e_{d-1})
\tn{d}s_1\hdots\tn{d}s_{d-1}\tn{d}t~.
\end{equation*}
We now prove \eqref{-15} by contradiction and assume that there exits $z^*\in(0,P_eX_{(N)}]$
such that
\beq{-3}
\int_{0}^{z^*} \tilde{f}(t)\tn{d}t  = \frac{G(z^*)}{G(P_eX_{(N)})}
>\frac{(z^*)^d}{(P_eX_{(N)})^d}=\frac{d}{(P_eX_{(N)})^{d}}\int_0^{z^*}t^{d-1}\tn{d}t,\eeq
where
$\tilde{f}(t):=\frac{ G^\prime (t) }{G(P_eX_{(N)})}$
 is the density of $\tilde F$.
 From \eqref{-3} and the monotonicity of the integral it follows that there exists   a point $t_0\in(0,z^*]$ with
\beq{-4}\tilde{f}(t_0)>\frac{d t_0^{d-1}}{(P_eX_{(N)})^d}.\eeq

For the following discussion, we introduce an  alternative  parametrization of the wedge $K$. As
  $\{{ e},{ e}_1,\hdots,{ e}_{d-1}\}$ is an  orthonormal basis of $\bb{R}^d$,
   every $x\in K$ can be represented as
$x=x_0+t_1{e}+s_1{ e}_1+\hdots+s_{d-1}{ e}_{d-1}$
for some $t_1 \in(0,l]$ and $s_i\in[-t_1 \tan(\vp),t_1 \tan(\vp)]$ for $i=1,\hdots,d-1$.
With the notation
\beq {ttilde}
{ e}_0:=\frac{t_1{ e}+s_1{ e}_1+\hdots+s_{d-1}{ e}_{d-1}}{\|t_1{ e}+s_1{ e}_1+\hdots+s_{d-1}{ e}_{d-1}\|}\quad\tn{and}
\quad\td{t}_1:=\|t_1{ e}+s_1{ e}_1+\hdots+s_{d-1}{ e}_{d-1}\|,
\eeq
we have $x=x_0+\td{t}_1{ e}_0$,
and the mapping of  $(t_1,s_1,\hdots,s_{d-1})$ to $(\td{t}_1,e_0)$ defines a bijection.
Hence,  any $x\in K$ can also be  uniquely represented by the vector ${ e}_0$ and the scalar $\td{t}_1$ (see Figure \ref{AA2}).
\input{exa117.tex}\vspace{-1cm}
Let $t_1\in [t_0, P_eX_{(N)}]$ and consider a point
\[
x=x_0+t_1e+s_1e_1+\hdots+s_{d-1}e_{d-1} = x_0+\td{t}_1e_0 \in K_N,
\]
where $e_0$ and $\td{t}_1$ are defined in \eqref{ttilde}.
Let $\td{t}_0=\frac{t_0}{t_1}\|t_1e+s_1e_1+\hdots+s_{d-1}e_{d-1}\|\leq\td{t}_1$ and define
\[y:=x_0+\td{t}_0e_0=x_0 +\td{t}_0\frac{t_1e+s_1e_1+\hdots+s_{d-1}e_{d-1}}{\|t_1e+s_1e_1+\hdots+s_{d-1}e_{d-1}\|}. \]
A straightforward calculation shows that $y=x_0+t_0e+\td{s}_1e_1+\hdots+\td{s}_{d-1}e_{d-1},$
where  $\td{s}_i  =\frac{t_0}{t_1}s_i,$ ($i=1,\hdots,d-1$).
Note that  $P_ey=t_0$ and that $f(y)\leq f(x)$, as $f$ is increasing on the wedge $K$. We therefore obtain
\bal \tilde{f}(t_1)&=\int_{- t_1\tan(\varphi)}^{ t_1\tan(\varphi)}\hdots\int_{-t_1\tan(\vp)}^{t_1\tan(\vp)}\frac{{f}(x_{0}+t_1e+s_1e_1+\hdots+s_{d-1}e_{d-1})}{G(P_eX_{(N)})}\tn{d}s_1\hdots\tn{d}s_{d-1}\\
&\geq \int_{- t_1\tan(\varphi)}^{ t_1\tan(\varphi)}\hdots\int_{-t_1\tan(\vp)}^{t_1\tan(\vp)}\frac{{f}(x_{0}+t_0e+
\frac{t_0}{t_1}s_1e_1+\hdots+\frac{t_0}{t_1}s_{d-1}e_{d-1})}{G(P_eX_{(N)})}\tn{d}s_1\hdots\tn{d}s_{d-1}. \\
& = \left(\frac{t_1}{t_0}\right)^{d-1}\int_{- t_0\tan(\varphi)}^{ t_0\tan(\varphi)}\hdots\int_{-t_0\tan(\vp)}^{t_0\tan(\vp)}\frac{{f}(x_{0}+t_0e+\td{s}_1e_1+\hdots+\td{s}_{d-1}e_{d-1})}{G(P_eX_{(N)})}\tn{d}\td{s}_1\hdots\tn{d}\td{s}_{d-1} \\
& = \left(\frac{t_1}{t_0}\right)^{d-1}\tilde{f}(t_0) .
\end{align*}
Using \eqref{-4}  we have
$
 \tilde{f}(t_1)  >
\big (\frac{t_1}{t_0}\big )^{d-1}\frac{d t_0^{d-1}}{(P_eX_{(N)})^d}=\frac{d t_1^{d-1}}{(P_eX_{(N)})^d}
$
for any point $t_1\in [t_0, P_eX_{(N)}]$, and
from \eqref{-3} we conclude
 \bal
 1=\int_{0}^{P_eX_{(N)}}\tilde{f}(t)\tn{d}t
 >\int_{0}^{z^*}\frac{dt^{d-1}}{(P_eX_{(N)})^d}\tn{d}t+\int_{z^*}^{P_eX_{(N)}}\frac{ dt^{d-1}}{(P_eX_{(N)})^d}\tn{d}t
 =\int_{0}^{P_eX_{(N)}}\frac{ dt^{d-1}}{(P_eX_{(N)})^d}\tn{d}t=1,\end{align*}
which is a contradiction
and proves the first  assertion  of Theorem \ref{01}.
\end{proof}

\begin{proof}  [\textbf{Proof of the representation  \eqref{P7}}] \relax
 It follows from Theorem \ref{P1} that for fixed $i$, the random variables  $\td{F}^i(P_{e_n^i}X_{1}^i),\ldots,\td{F}^i(P_{e_n^i}X^i_{N^i-1})$ are independent and
 uniformly distributed on the interval $[0,1],$  given $\{X_j^i\in K^i_{N^i}\} \ (j=1,\ldots,N^i-1)$ and $N^i$. Here,
  $\td{F}^i$ denotes the conditional distribution function of $P_{e_n^i}X$ on $K^i$.
 Recall that the wedges $K_{N_i}^i,\,i=1,\ldots,M_n,$ are disjoint. Standard arguments and   the law of iterated expectations  yield stochastic independence of observations lying in different wedges.
\end{proof}

\begin{proof}[\textbf{Proof of Theorem \ref{P-8}}]\relax
Let $\bs{a}\in\{-1,1,0\}^{M_n}$ be the vector determined by
\[a_i=\begin{cases} 1,&\tn{if }f\tn{ is increasing on }K^i,\\
-1,&\tn{if }f\tn{ is decreasing on }K^i,\\
0,&\tn{ else.}\end{cases}\]
 The probability of at least one false rejection among
all tests in \eqref{incrrej} and \eqref{decrrej} can be estimated by Theorem  \ref{01}, that is
\begin{align*}
p& := ~ \bb{P}_f\big(\exists i\in\{i\,|\,a_i\neq 0\}:a_iT_{K^i}<-\td{c}_{K^i}(\alpha)\big|N^1,\hdots,N^{M_n}\big)\\
& \leq ~ \bb{P}_f\big(\exists i\in\{i\,|\,a_i\neq 0\}:a_iT_{K^i}^U<-\td{c}_{K^i}(\alpha)\big|N^1,\hdots,N^{M_n}\big)\\
& \leq ~1-\bb{P}\left(|T_{K^i}^U|\leq  \td{c}_{K^i}(\alpha)\tn{ for all }i=1,\hdots,M_n\big|N^1,\hdots,N^{M_n}\right).
\end{align*}
Using \eqref{m1} we further deduce
\begin{align*}
p \leq&~Ê1-\bb{P}\Big(\Big(\sqrt{\frac{3}{N^i-1}}\left|T_{K^i}^U\right|-\Gamma\Big (\frac{N^i}{n-1}\Big)\Big)\leq
\td{\kappa}_n(\alpha)\tn{ for all }i=1,\hdots,M_n\big|N^1,\hdots,N^{M_n}\Big )\\
=&~1-\bb{P}\Big(\max_{i=1,\hdots,M_n}\Big(\sqrt{\frac{3}{N^i-1}}\left|T_{K^i}^U\right|-\Gamma\Big (\frac{N^i}{n-1}\Big)\Big)
\leq \td{\kappa}_n(\alpha)\big|N^1,\hdots,N^{M_n}\Big ) ~=~\alpha.\end{align*}
\end{proof}

\subsection{Further Results}
This section provides a general  consistency result which is the main ingredient for the proof of Theorems \ref{S1A} and \ref{S1B}.
The consistency result   stated in Theorem \ref{98} below is a more general result.
The following notation   is used throughout this section for the precise statement of the assumptions on the density $f$.

\begin{defi} \label{d61}
Let $K$ be a wedge with vertex $x_0$ and $j\in\{1,2\}$. We define
\begin{eqnarray}  \label{supf}
\inf_K f^{(j)}&:=&\inf_{x_{0}+\td{t}e_0\in K, 0\leq \td{s}<\td{t}  }\frac{f(x_{0}+\td{t}e_0)-f(x_{0}
+\td{s}e_0)}{\td{t}^j-\td{s}^j}, \\ \nonumber
\sup_K f^{(j)}&:=&\sup_{x_{0}+\td{t}e_0\in K, 0\leq \td{s}<\td{t}  }\frac{f(x_{0}+\td{t}e_0)-f(x_{0}
+\td{s}e_0)}{\td{t}^j-\td{s}^j}.
\end{eqnarray}
\end{defi}
For a better interpretation of  Definition \ref{d61}, let $f$ be twice continuously differentiable in a neighbourhood of
$x_0$. A straightforward application of   the mean value theorem shows
\[\inf_K f^{(1)}=\inf_{x_{0}+\td{t}e_0\in K }\langle\nabla f(x_0+\td{t}e_0),e_0\rangle,
\quad\sup_K f^{(1)}=\sup_{x_{0}+\td{t}e_0\in K }\langle\nabla f(x_0+\td{t}e_0),e_0\rangle.\]
If we have a mode in $x_0$,   the gradient  $\nabla f(x_0)$ vanishes and  $\sup_K f^{(1)}$
vanishes as well. Similarly,   a Taylor expansion of order 2 yields
\[\inf_K f^{(2)}=\frac{1}{2}\inf \big \{ \langle e_0,H_f(x_0)e_0\rangle \big | \{x_{0}+\td{t}e_0,\td{t}\geq0\}\cap K\neq\emptyset \} \big | \big \} +o(1)\]
\[\sup_K f^{(2)}=\frac{1}{2} \sup \big \{ \langle e_0,H_f(x_0)e_0\rangle \big | \{x_{0}+\td{t}e_0,\td{t}\geq0\}\cap K\neq\emptyset \} \big | \big \} +o(1)\]
for $l\rightarrow0$, where $H_f(x_0)$ denotes the Hessian of $f$ in $x_0$.
In the situation of Theorem \ref{S1A}, the condition $e_0^\top H_f(x_0)e_0\leq -c<0$ holds for all $e_0\in\bb{R}^d$ with
$\|e_0\|=1$. Thus, only the case  $j=2$ is relevant for its proof. However,
if the assumption is violated and the Hessian is indefinite, the results can be generalized in a straightforward
manner by considering orders $j\geq 3$ as well. For the proof of Theorem \ref{S1B}, the case $j=1$ will be used.


To simplify notation, let
\[\td{\vp}=\tan(\vp)\tn{ for }\vp\in\big(0,\frac{\pi}{2}\big)\quad\tn{ and }\quad F(K)=\int_{K}f (\bs{x})\tn{d}\bs{x}.\]
The  quantities
\begin{eqnarray*}
&& H_+^j(f,{K}):=\frac{(2\td{\varphi})^{d-1}l^{d+j}\inf_{{K}}f^{(j)}}{\sqrt{F(K)}}~,~~
 {H}_-^j(f,{K}):=\frac{(2\td{\varphi})^{d-1}l^{d+j}\sup_{{K}}f^{(j)}}{\sqrt{F(K)}}
 \quad (j=1,2)
 \end{eqnarray*}
depend  on the size of the wedge $K$ through its length $l$ and angle $\vp$ and on the monotonicity
of $f$ on $K$  and   are the key objects in the following discussion. We begin by showing
that the quantiles $\td{\kappa}_n(\alpha)$
and   $\ol{\kappa}_n(\alpha)$ defined in    \eqref{m1} and \eqref{m2} are bounded from above by a constant independent
of $M_n,N^i,i=1,\hdots,M_n,$ and $n$. As a consequence, the same holds for the quantiles $\td{\kappa}_n'(\alpha)$
defined in \eqref{glei1}.

\begin{satz}\label{C1}
There exists a constant $A > 0$, such that
$\max \{ \td{\kappa}_n(\alpha),\ol{\kappa}_n(\alpha)\} \leq A. $
\end{satz}

\begin{proof}[\textbf{Proof}]
We only consider $\td{\kappa}_n$ as the  result for $\ol{\kappa}_n$ can be shown similarly.
From the discussion in  Section \ref{2.2}, it follows  that
\[T_{K^i}^U\overset{d}{=}\sum_{j=1}^{N^i-1}\beta(U_j^i) \tn{ for }i=1,\hdots,M_n\]
(conditionally on  $N^1,\hdots,N^{M_n}$),  where $U_j^i,\;j=1,\hdots,N^i-1,i=1,\hdots,M_n,$ are  independent   
uniformly distributed  random variables
on the interval $[0,1]$.
Recall the definition of $T_{jk}(\bs{U})$ in \eqref{n1}, then we will show a the end of the proof that,
conditionally on
 $N^1,\hdots,N^{M_n}$,
\beq{-80}T_{\sum_{k=1}^{i-1}N^k,\sum_{k=1}^iN^k }(\bs{U})
= \sum_{\ell =\sum_{k=1}^{i-1}N^k+1}^{ \sum_{k=1}^iN^k-1}
 \beta \big ( U_{(\ell ; \sum_{k=1}^{i-1}N^k,\sum_{k=1}^iN^k) }\big )
\overset{d}{=}\sum_{j=1}^{N^i-1}\beta(U_j^i)
\eeq
($i=1,\hdots,M_n$).  As the statistic  $T_n(\bs{U})$  defined in \eqref{TnU} calculates the maximum over more scales
than the statistic defined in \eqref{m1}, we obtain  $\td{\kappa}_n \leq \kappa_n $,
where $\kappa_n $ is the $(1-\alpha)$-quantile of the statistic $T_n(\bs{U})$  defined in \eqref{TnU}.
By Theorem 3.1 in \cite{MR2435455}, there exists an upper bound $A>0$ for  $\kappa_n$, which is  independent of
$n$. This completes the proof.\\ 
For a proof of   \eqref{-80}, we finally note that for  a sample of independent uniformly distributed random variables
$U_1, \ldots, U_n$ on the interval $[0,1]$ and  fixed $1\leq j<k\leq n$ with   $k-j>1$, the random variables
$
U_{(l;j,k)} ~=~ ( {U_{(l)}-U_{(j)}})/({U_{(k)}-U_{(j)}}) $
($l=j+1,\hdots,k-1$) have the same distribution as the order statistic of $k-j-1$ independent
 uniformly   distributed random variables on the interval $[0,1]$. 
 \end{proof}

Following the notation from Section  \ref{2.2} (recall that $
\Gamma(\delta)=\sqrt{2\log ( \tfrac{\exp(1)}{\delta}  )} $), we define the quantiles
\[{c}_{K_n^i}(\alpha)=\sqrt{\frac{N^i-1}{3}}\Big(A+\Gamma\Big(\frac{N^i}{n-1}\Big)\Big)\] $( i=1,\hdots,M_n)$
and provide   a general consistency result   for locally increasing and
decreasing densities.

\begin{satz}\label{98}
Let $j\in\{1,2\}$ and  $\ca{K}_n = \{K^1_n,\ldots,K_n^{M_n}\}$   be a family of $M_n$ pairwise disjoint wedges with length $l_n>0$ and angle
$\vp_n\in(0,\frac{\pi}{2})$.
\begin{enumerate}
\item If (for all $i=1,\ldots,M_n$) the condition $
H_+^j(f,{K^i_n})\geq D\big (\frac{\Gamma(F(K_n^i))}{\sqrt{2}}+b_n\big){\frac{1}{\sqrt{n}}}
$
holds for some  constant
\begin{equation}\label{63i}
D>\frac{j(d+j)2\sqrt{2}}{{(2d+j)\Big (-1+\big\{{1+\frac{2j^2}{(2d+j)^2}}\big\}^\frac{1}{2}\Big)\Big({1-\frac{j^2}{2(2d+j)^2}
\big\{{-1+\big[{1+4\frac{(2d+j)^2}{j^2}}}\big]^\frac{1}{2}\big\}}\Big)^\frac{1}{2}}}
\end{equation}
and a positive sequence $b_n$ satisfying    $b_n\rightarrow\infty$ and $b_n=o(\sqrt{\log(n)})$
as  $n\rightarrow\infty$, then
\[\bb{P}\big(T_{K^i_n}>c_{K^i_n}(\alpha)\text{ for all } i=1,\ldots,M_n \big)\rightarrow 1.\]
\item If (for all $i=1\ldots,M_n$) the condition
 $ {H}_-^j(f,{K^i_n})\leq- D\big (\frac{\Gamma(F(K_n^i))}{\sqrt{2}}+b_n\big )\frac{1}{\sqrt{n}}$
 holds for some constant
\[D>\frac{2\sqrt{2}(2d+j)(d+j)}
{{j\big(1-\frac{d}{d+j}\big)^{(d+j)/d}
\Big[{1-\frac{d^2}{2(2d+j)^2}{\big(-1+\big\{{1+4(\frac{2d+j}{d})^2}\big\}^{1/2}\big)\Big]^{1/2}}}}}\]

and a positive sequence $b_n$ satisfying    $b_n\rightarrow\infty$ and $b_n=o(\sqrt{\log(n)})$
as  $n\rightarrow\infty$, then
\[\bb{P}\big(T_{K_n}<-c_{K_n}(\alpha)\text{ for all }{i=1,\ldots,M_n} \big)\rightarrow 1.\]
\end{enumerate}
\end{satz}
\begin{bem}{\rm
It follows from the proof of Theorem \ref{98} below, that the bounds on $H_+^j(f,{K_n})$ resp. $H_-^j(f,{K_n})$
imply in particular that $M_n=o(n)$ as  $n\rightarrow\infty$.}
\end{bem}
The proof of Theorem \ref{98} is divided into eight parts: seven technical lemmas (Lemma \ref{02} - Lemma \ref{05})
and the main part of the proof. We first state and prove the technical lemmas and finally combine the results in order 
to  complete the proof of  Theorem \ref{98}.
 For each of the lemmas, we assume  that the conditions of Theorem \ref{98} hold. At first, we consider only one  fixed  
 wedge $K$ with   length $l_n=l$ and angle $\vp_n=\vp$. To simplify notation, let
\bal
\ol{F}^j_\pm(z)&:=z+c_\pm \left(\frac{2^{d-1}}{d+j} z^\frac{d+j}{d}-\frac{2^{d-1}}{d+j}z\right)
\end{align*}
for $j\in\{1,2\}$ and $z\in[0,1]$ with
\[c^j_+:=\frac{(\tan(\varphi))^{d-1}(P_eX_{(N)})^{d+j}\inf_K f^{(j)}}{F(K_N)}\quad\tn{and}
\quad{c}^j_-:=\frac{(\tan(\varphi))^{d-1}(P_eX_{(N)})^{d+j}\sup_K f^{(j)}}{F(K_N)}.\]

\begin{lemma}\label{02}  ~
\begin{enumerate}
\item If $\inf_K f^{(j)}>0$, then  $ \tilde{F}(z)\leq \ol{F}^j_+\big(\frac{z^d}{(P_eX_{(N)})^d}\big) $
for $z\in(0,P_eX_{(N)}]$.
\item If $\sup_K f^{(j)}<0$, then $\tilde{F}(z)\geq \ol{F}^j_-\big(\frac{z^d}{(P_eX_{(N)})^d}\big)$
for $z\in(0,P_eX_{(N)}]$.
\end{enumerate}
\end{lemma}

\begin{proof}[\textbf{Proof}]
We only prove the first part   and define
 the auxiliary function
\bal\overline{f}^j_+(s_1,\hdots,s_{d-1},t):=&\mathbbm{1}_{K_N} (x_0+te+s_1e_1+\hdots+s_{d-1}e_{d-1} ) \\
&\quad\cdot\Big(\frac{1}{|K_N|}+\Big(\frac{\inf_Kf^{(j)}}{F(K_N)}\Big)\Big(t^j-\frac{d}{d+j}(P_eX_{(N)})^j\Big)\Big),\end{align*}
where $|K_N|:=\int_{K_N}1\tn{d}\bs{x}$ denotes the volume of $K_N$. Note that
\begin{align}\begin{split}\label{-61}&\int_0^{z}\int_{-t\tan(\varphi)}^{t\tan(\varphi)}\hdots
\int_{-t\tan(\vp)}^{t\tan(\vp)}\ol{f}^j_+(s_1,\hdots,s_{d-1},t)\tn{d}s_1\hdots\tn{d}s_{d-1}\tn{d}t\\
=&\frac{z^d}{(P_eX_{(N)})^d}+\int_0^{z}(2\tan(\varphi)t)^{d-1}\Big(\frac{\inf_Kf^{(j)}}{F(K_N)}\Big)\Big(t^j
-\frac{d}{d+j}(P_eX_{(N)})^j\Big)\tn{d}t
\\=&\frac{z^d}{(P_eX_{(N)})^d}+c_+\Big(\frac{2^{d-1}}{d+j}\frac{z^{d+j}}{(P_eX_{(N)})^{d+j}}-\frac{2^{d-1}}{d+j}
\frac{z^d}{(P_eX_{(N)})^d}\Big) = \overline{F}^j_+\Big(\frac{z^d}{(P_eX_{(N)})^d}\Big).\end{split}\end{align}
In particular, as $\ol{F}^j_+(1)=1$, the function $\overline{f}^j_+$ defines a density on $K_N$.
We now prove  assertion \textit{(i)} by
 contradiction and
assume that there exists $z^*\in(0,P_eX_{(N)}]$,  such that
\beq{-9}\tilde F (z^*)= \int_0^{z^*} \td{f}(t)dt =\frac{G(z^*)}{G(P_eX_{(N)})}>\overline{F}^j_+\Big(\frac{(z^*)^d}{(P_eX_{(N)})^d}\Big) = \int_0^{z^*}h^j_+(t)dt, \eeq
where $\td {f}$ and
\[h^j_+(t):=\int_{-t\tan(\varphi)}^{t\tan(\varphi)}\hdots\int_{-t\tan(\vp)}^{t\tan(\vp)}\ol{f}^j_+(s_1,\hdots,s_{d-1},t)\tn{d}s_1\hdots\tn{d}s_{d-1}\]
denote the density of $\td{F}$ and $\ol{F}^j_+$, respectively.
Due to the monotonicity of the integral, there exists  a point $t_0\in(0,z^*]$ with
\beq{l111} \td{f}(t_0)>h_+^j(t_0) ,\eeq
  which implies
\begin{align}\begin{split}\label{11}\td{f}(t_1)>h^j_+(t_1)\tn{ for all } t_1\in[t_0, P_eX_{(N)}].\end{split}\end{align}
 For a proof of \eqref{11},  let $t_1>t_0$ and
$x=x_0+t_1e+s_1e_1+\hdots+s_{d-1}e_{d-1}\in K_N$.
As in  the proof of Theorem \ref{01}, we use the representation
$x=x_0+\td{t}_1e_0$,
where $\td{t}_1$ and $e_0$ are defined in \eqref{ttilde}.
For $\td{t}_0:=\frac{t_0}{t_1}\td{t}_1<\td{t}_1$, let
\[y:=x_0+\td{t}_0e_0=x_0+t_0e_0+\frac{t_0}{t_1}s_1e_1+\hdots+\frac{t_0}{t_1}s_{d-1}e_{d-1}\in K_N.\]
Using $\inf_Kf^{(j)}>0$, we find
\begin{eqnarray}\label{-12}
f(x)-f(y)&=& f(x_0+\td{t}_1e_0)-f(x_0+\td{t}_0e_0)
=(\td{t}_1^j-\td{t}_0^j)\frac{f(x_0+\td{t}_1e_0)-f(x_0+\td{t}_0e_0)}{\td{t}_1^j-\td{t}_0^j}
\\ \nonumber
&\geq& (\td{t}_1^j-\td{t}_0^j)\inf_Kf^{(j)}=\td{t}_1^j\Big(1-\frac{t_0^j}{t_1^j}\Big)\inf_Kf^{(j)}
\geq (t_1^j-t_0^j)\inf_Kf^{(j)},
\end{eqnarray}
where the last estimate follows since $\td{t}_1\geq t_1.$
Recall that $G(P_eX_{(N)})=F(K_N)$, then we obtain
\begin{align*}
 \td{f}(t_1)&=\int_{-t_1\tan(\varphi)}^{t_1\tan(\varphi)}\hdots
\int_{-t_1\tan(\vp)}^{t_1\tan(\vp)}\frac{f(x_{0}+t_1e+s_1e_1+\hdots +s_{d-1}e_{d-1})}{G(P_eX_{(N)})}\tn{d}s_1
\hdots\tn{d}s_{d-1}\\&\geq\int_{-t_1\tan(\varphi)}^{t_1\tan(\varphi)}\hdots\int_{-t_1\tan(\vp)}^{t_1\tan(\vp)}
\frac{1}{F(K_N)}\Big( f(x_{0}+t_0e+\textstyle{\frac{t_0}{t_1}}s_1e_1+\hdots+\frac{t_0}{t_1}s_{d-1}e_{d-1})\\
&\hspace{6.3cm}+(t_1^j-t_0^j)\inf_Kf^{(j)}\Big)\tn{d}s_1\hdots\tn{d}s_{d-1}.\end{align*}
A change of variables yields
\begin{align*} \td{f}(t_1)\geq & \left(\frac{t_1}{t_0}\right)^{d-1}
\int_{-t_0\tan(\varphi)}^{t_0\tan(\varphi)}\hdots
\int_{-t_0\tan(\vp)}^{t_0\tan(\vp)}\frac{ f(x_{0}+t_0e+\tilde{s}_1e_1+\hdots+\tilde{s}_{d-1}e_{d-1})}{F(K_N)}
\tn{d}\tilde{s}_1\hdots\tn{d}\td{s}_{d-1}\\&+(2\tan(\varphi)t_1)^{d-1}(t_1^j-t_0^j)\frac{\inf_Kf^{(j)}}{F(K_N)}\\
=&\left(\frac{t_1}{t_0}\right)^{d-1}\td{f}(t_0)+(2\tan(\varphi)t_1)^{d-1}(t_1^j-t_0^j)\frac{\inf_Kf^{(j)}}{F(K_N)},
\end{align*}
and straightforward calculations show that
  \begin{align*}
    h_+^j(t_1) & =\left(\frac{t_1}{t_0}\right)^{d-1}h_+^j(t_0)+(2\tan(\varphi)t_1)^{d-1}(t_1^j-t_0^j)
    \frac{\inf_Kf^{(j)}}{F(K_N) } \\
& <\left(\frac{t_1}{t_0}\right)^{d-1}\td{f}(t_0)
+(2\tan(\varphi)t_1)^{d-1}(t_1^j-t_0^j)\frac{\inf_Kf^{(j)}}{F(K_N)}\leq\tilde f(t_1),
\end{align*}
where we used \eqref{l111} to obtain the strict inequality.
 From \eqref{-9} and \eqref{11}  we  also get
$$
\int_{0}^{P_eX_{(N)}}\td{f}(t)\tn{d}t   >\int_{0}^{z^*}h^j_+(t)\tn{d}t+\int_{z^*}^{P_eX_{(N)}}h^j_+(t)\tn{d}t  =\ol{F}^j_+(1)=1,
 $$
which contradicts the condition $1=\td{F}(P_eX_{(N)})=\int_{0}^{P_eX_{(N)}}\td{f}(t)\tn{d}t$.
This completes the proof of Lemma \ref{02}.
\end{proof}

\begin{lemma}\label{1109} ~
\begin{enumerate}
\item If $\inf_K f^{(j)}>0$, we have $ \sum_{i=1}^{N-1} \beta\big((\ol{F}_+^j)^{-1}(\tilde{F}(P_eX_{(i)}))\big)\leq T_{K}.$
\item If $\sup_K f^{(j)}<0$, we have $  \sum_{i=1}^{N-1} \beta\big((\ol{F}_-^j)^{-1}(\tilde{F}(P_eX_{(i)}))\big)\geq T_{K}. $
\end{enumerate}
\end{lemma}

 \begin{proof}[\textbf{Proof}]
 We only prove the first  part and begin showing that the function   $z\mapsto \ol{F}_+^j(z)$ is strictly increasing
for $z\in[0,1]$. Recalling the representation
\eqref{-61},
it is sufficient to prove that  the inequality
\[\ol{f}_+^j(s_1,\hdots,s_{d-1},t)> 0\quad\tn{for }t\in(0,P_eX_{(N)}]  \]
holds for all $(s_1,\ldots,s_{d-1}) \in [-t \ \rm{tan}\varphi, t  \ \rm{tan} \varphi  ]$.
For the sake of simplicity, we suppress the dependence of $\ol{f}^j_+$ on $(s_1,\ldots,s_{d-1})$ and note that the function $t\mapsto\ol{f}_+^j(t) = \ol{f}^j_+(t,s_1,\ldots,s_{d-1})$
is strictly increasing. Therefore, it remains to show $\ol{f}^j_+(0)\geq 0$. We prove this inequality by contradiction and
assume that  $\ol{f}^j_+(0)<0$. For
$x=x_0+te+s_1e_1+\hdots+s_{d-1}e_{d-1}\in K_N$,
it follows (using $\td{t}_0=0$ in \eqref{-12})
\begin{eqnarray} \label{-13}
&& \frac{f(x)} {F(K_N)}
  \geq\frac{f(x)-f(x_0)}{F(K_N)}
 \geq t^j \frac{\inf_Kf^{(j)}}{F(K_N)}  =\ol{f}^j_+(t)-\ol{f}^j_+(0)
 > \ol{f}^j_+(t).
\end{eqnarray}
Integrating both sides of \eqref{-13}   leads to a contradiction.
Consequently, the map $z\mapsto \bar{F}^j_+(z)$ is strictly increasing on $[0,1]$,
which implies (using Lemma \ref{02} and the monotonicity of the function $\beta$) that
$\beta((\ol{F}^j_+)^{-1} (\td{F}(P_e X_{(j)}))) \leq \beta \big(\frac{(P_e X_{(j)})^d}{P_eX_{(N)})^d}\big)$, whenever
$\td{F}(P_eX_{(j)})\neq  0  $ and $P_eX_{(j)}\neq P_eX_{(N)}$. However, it is easy to see that these cases correspond to
$P_eX_{(j)}=0$ and $\td{F}(P_eX_{(j)})=1$, where there is in fact equality. Thus the proof of the first part is completed.  
\end{proof}

The conditional expectation considered in the following lemma is used to derive a bound on the corresponding  conditional probability via the Hoeffding inequality in Lemma \ref{05} below.

\begin{lemma}\label{03} ~
\[ \leqno{(i)} ~~~~~~~~~~~~~~
\bb{E}\Big(\sum_{i=1}^{N-1} \beta\big((\ol{F}_+^j)^{-1}(\tilde{F}(P_eX_{(i)}))\big)|N,X_{(N)}\Big)
=\frac{2^{d-1}j(N-1)}{(2d+j)(d+j)}c_+^j.\]
\[\leqno{(ii)} ~~~~~~~~~~~~~~
\bb{E}\Big(\sum_{i=1}^{N-1} \beta\big((\ol{F}_-^j)^{-1}(\tilde{F}(P_eX_{(i)}))\big)|N,X_{(N)}\Big)
=\frac{2^{d-1}j(N-1)}{(2d+j)(d+j)}\td{c}_-^j.\]

\end{lemma}

\begin{proof}[\textbf{Proof}]
We only prove the first part. Let $U_1,\hdots,U_{N-1}$ be independent uniformly distributed random variables on the interval $[0,1]$. Theorem \ref{P1} yields
\[\sum_{i=1}^{N-1} \beta\left((\ol{F}_+^j)^{-1}(\tilde{F}(P_eX_{(i)}))\right)\overset{d}{=}
\sum_{i=1}^{N-1} \beta\left((\ol{F}_+^j)^{-1}(U_i)\right), \]
given $N$ and $X_{(N)}$.
The assertion now follows
from
\bal \bb{E}\big(\beta((\ol{F}_+^j)^{-1}(U_1))|N,X_{(N)}\big)&=\int_0^1 \beta(x)g_+^j(x)\tn{d}x
 =\frac{2^{d-1}j}{(2d+j)(d+j)}c_+^j,\end{align*}
 where $  g_+^j(x):=\frac{\tn{d}}{\tn{d}x}\ol{F}^j_+(x)=1+c_+^j\big(\frac{2^{d-1}}{d}x^{\frac{d+j}{d}-1}-\frac{2^{d-1}}{d+j}\big)$.
\end{proof}

In the following, we consider a sequence of wedges $(K_n)_{n\in\bb{N}}$  given by vertices $x_0^n$,
directions $e^n$, lengths $l_n>0$ and angles $\vp_n\in(0,\frac{\pi}{2})$.
Furthermore, we denote by $e_1^n,\hdots,e_{d-1}^n$ the orthonormal basis of $(\mathrm{span}\{e^n\})^\perp$
and let $\delta_n=F(K_n),\;\delta_{N_n}=F(K_{N_n})$.
Lemma \ref{04} and Lemma \ref{07} below ensure the feasibility of our procedure in an asymptotic sense. They show that the random
wedge $K_{N_n}$ is of similar size than the deterministic, predefined wedge $K_n$, i.e. that  its complement  $K_n\backslash K_{N_n}$ is small.
Note that the test procedure can only be consistent  if the wedge $K_n$ contains a sufficiently large number of
observations. Therefore, we introduce for $\gamma \in (0,\frac{1}{2}]$ and  $0<\ve<1$ the conditional probability ${\bb{P}}_{\mathcal{N}_n}$  given  the event
\begin{align}\label{Nn}
\mathcal{N}_n={\Big\{N_n\geq(1-\gamma)n\delta_n,\;\frac{|K_{N_n}|}{|{K_n}|}\geq1-\ve\Big\}}.
\end{align}
The results of Lemma \ref{05} below are only shown for ${\bb{P}}_{\mathcal{N}_n}$.
However, the following Lemmas \ref{04} and \ref{07} demonstrate that these conditions are asymptotically negligible.
For example,   Lemma \ref{04} shows that, with increasing $n$, the wedge $K_{N_n}\subseteq K_n$ approximates $K_n$ in probability
at an exponential rate.

 \begin{lemma}\label{04}
Let  $j\in\{1,2\}$ and $0<\ve<1$.
\begin{enumerate}
 \item If the assumptions of Theorem \ref{98} \textit{(i)} are satisfied, then 
 \[
 \bb{P}\Big(N_n=0\tn{ or }\frac{|K_{N_n}|}{|{K_n}|}<1-\ve\Big)
\leq
 \exp\Big(-\frac{D}{d+j}
\sqrt{\frac{n\delta_n}{2}}\Gamma(\delta_n)\big (1-(1-\ve)^\frac{d+j}{d}\big)\Big).\]
\item  If the assumptions of Theorem \ref{98} \textit{(ii)} are satisfied, then
\[
\bb{P}\Big(N_n=0\tn{ or }\frac{|K_{N_n}|}{|{K_n}|}<1-\ve\Big)
\leq \exp\Big(-\frac{D}{d+j}
\sqrt{\frac{n\delta_n}{2}}\Gamma(\delta_n)\big ((1-\ve)^\frac{d+j}{d}+\frac{d+j}{d}\ve-1 \big )\Big).\]
\end{enumerate}
\end{lemma}

\begin{proof}[\textbf{Proof}]  
We only prove the first part. Since
$\frac{|K_{N_n}|}{|{K_n}|}=\frac{(P_{e^n}X_{(N_n)})^d}{l_n^d}< 1-\ve$
 if and only if
$P_{e^n}X_{(N_n)}< l_n(1-\ve)^\frac{1}{d}$,
we obtain
\begin{align*}
\Big\{N_n=0\tn{ or }\frac{|K_{N_n}|}{|{K_n}|}<1-\ve\Big\}=
\big\{N_n=0\tn{ or }P_{e^n}X_{(N_n)}< l_n(1-\ve)^\frac{1}{d}\big\}.
\end{align*}

Define
 \[I_{r_n}:=K_n\cap\left\{x\in\bb{R}^d:\;P_{e^n}x\geq l_n(1-\ve)^\frac{1}{d}\right\}  \]
 (see Figure \ref{abbi11}), then
$
\{N_n=0\tn{ or }\frac{|K_{N_n}|}{|{K_n}|}<1-\ve\big\}\subseteq\{\tn{no observation in }I_{r_n}\}.$
\begin{figure}
\input{irn.tex}
\caption{\emph{The section $I_{r_n}.$}}
\label{abbi11}
\end{figure}
Now, recall for a proof of part \textit{(i)}   from \eqref{-12}  that (with $t_0=0$ and $t_1=t$) $f(x_{0}^n+te^n+s_1e^n_1+\hdots+s_{d-1}e^n_{d-1})
\geq t^j\inf_{K_n}f^{(j)}$.  Thus,
  \begin{align*} p_{f,I_{r_{n}}}:&=\int_{I_{r_n}}f(\bs{x})\tn{d}\bs{x}\\
  &=\int_{l_n(1-\ve)^\frac{1}{d}}^{l_n}\int_{-t\td{\vp}_n}^{t\td{\vp}_n}\hdots\int_{-t\td{\vp}_n}^{t\td{\vp}_n}f(x_{0}^n+te^n+s_1e^n_1+\hdots+s_{d-1}e^n_{d-1})\tn{d}s_1\hdots\tn{d}s_{d-1}\tn{d}t\\
 &\geq\int_{l_n(1-\ve)^\frac{1}{d}}^{l_n}\int_{-t\td{\vp}_n}^{t\td{\vp}_n}\hdots\int_{-t\td{\vp}_n}^{t\td{\vp}_n}t^j
  \inf_{K_n}f^{(j)}\tn{d}s_1\hdots\tn{d}s_{d-1}\tn{d}t
 = \frac{(2\td{\vp}_n)^{d-1}}{d+j}l_n^{d+j}\big(1-(1-\ve)^\frac{d+j}{d}\big)\inf_{{K_n}}f^{(j)}.\end{align*}
The assumptions of  Theorem \ref{98} \textit{(i)} imply that  $H_+^j(f,{K_n})\geq D\frac{\Gamma(\delta_n)}{\sqrt{2n}}
$
and therefore
$ p_{f,I_{r_{n}}}\geq\frac{D}{d+j}\sqrt{\frac{\delta_n}{2n}}\Gamma(\delta_n) (1-(1-\ve)^\frac{d+j}{d}).$
As the variables
$Z_i= \mathbbm{1}_{I_{r_n}}(X_i)$ are Bernoulli distributed with parameters $p_{f,I_{r_{n}}}$, we have
\bal  \bb{P}(\tn{no observation in }I_{r_n})
 \leq \exp\Big(-\tfrac{D}{d+j}\sqrt{\tfrac{n\delta_n}{2}}\Gamma(\delta_n)\big(1-(1-\ve)^\frac{d+j}{d}\big)\Big).\end{align*}
\end{proof}

As the number  $N_n$ of observations in $K_n$ is Bin$(n,\delta_n)$-distributed, we obtain the following result from Chernoff's Inequality.
\begin{lemma}\label{07}
$ {\bb{P}}\left( N_n\leq(1-\gamma)n\delta_n\right)\leq\exp\big(-n\delta_n\frac{\gamma^2}{2}\big) $
for any $\gamma \in (0,\frac{1}{2}]$.
\end{lemma}

\begin{lemma}\label{06} Let $j\in\{1,2\}$.
\begin{enumerate}
\item If the assumptions of Theorem \ref{98} \textit{(i)} are satisfied, then
\begin{enumerate}
\item $\quad n\delta_n\geq\displaystyle{\frac{D^2}{(d+j)^2}\frac{\Gamma(\delta_n)^2}{2}=:\td{c}^j_+\frac{\Gamma(\delta_n)^2}{2}}$,
\item $\quad n\delta_n\geq L_n\td{c}^j_+\log(\exp(1)n),\tn{ where }L_n\geq1-o(1)$ for $n\rightarrow\infty$.

\end{enumerate}
 \item If the assumptions of Theorem \ref{98} \textit{(ii)} are satisfied, then
\begin{enumerate}
\item[\textit{(a')}]$\quad n\delta_n\geq\displaystyle{\frac{j^2D^2}{(d(d+j))^2}\frac{\Gamma(\delta_n)^2}{2}
=:\td{c}^j_-\frac{\Gamma(\delta_n)^2}{2}}$.
\item[\textit{(b')}] $\quad n\delta_n\geq L_n\td{c}^j_-\log(\exp(1)n),\tn{ where }L_n\geq1-o(1)$ for $n\rightarrow\infty$.
\end{enumerate}
\end{enumerate}
\end{lemma}

\begin{proof}[\textbf{Proof}]
We only prove the first part. As in the proof of Lemma \ref{04}, we obtain
\[\delta_n=\int_{K_n}f (\bs{x}) \tn{d}\bs{x}\geq\frac{(2\td{\vp}_n)^{d-1}}{d+j}l_n^{d+j}\inf_{{K_n}}f^{(j)}.\]
Hence,
\[H^j_+(f,{K_n})=\frac{(2\td{\vp}_n)^{d-1}l_n^{d+j}\inf_{{K_n}}f^{(j)}}{\sqrt{\delta_n}}\leq
\frac{(2\td{\vp}_n)^{d-1}l_n^{d+j}\inf_{{K_n}}f^{(j)}}{\frac{(2\td{\vp}_n)^{d-1}}{d+j}l_n^{d+j}\inf_{{K_n}}f^{(j)}}
\sqrt{\delta_n}=(d+j)\sqrt{\delta_n}.\]
Therefore, it follows from the assumption  $H^j_+(f,{K_n})$  that
$ n\delta_n\geq\tfrac{D^2}{(d+j)^2}\tfrac{\Gamma(\delta_n)^2}{2}.  $
Part \textit{(b)} is a consequence of Lemma 7.5 in \cite{MR2435455}.
\end{proof}

\bigskip

The following Lemma provides the key to prove  consistency.
Note that for the construction of the test statistic $T_{K_n}$,  it is necessary that at least two
observations are contained in the wedge $K_n$.
Given the event $\mathcal{N}_n$, we have that $N_n\geq 2$ if $n\delta_n\geq 4.$
 If the assumptions of Theorem \ref{98}
 \textit{(i)} hold, it follows from  \textit{(a)} that $n\delta_n\geq 4$ is fulfilled for $D\geq{2(d+j)}{}$.
Similarly, if  the assumptions of Theorem \ref{98} \textit{(ii)} hold, then \textit{(a')} yields the condition $D\geq \frac{2d(d+j)}{j}$.
\begin{lemma}\label{05}
Let $j\in\{1,2\}$ and ${\bb{P}}_{\mathcal{N}_n}$ denote the probability conditional on the event $\mathcal{N}_n$  defined in \eqref{Nn},
and define for
 $\gamma \in (0,\frac{1}{2}]$,  $0<\ve<1$ and $\eta>0$ the constant
\begin{align}\label{D}
\mathbf{D}^j(\eta,\delta_n):= \frac{ \frac{(2d+j)(d+j)}{j\sqrt{3}}}{(1-\ve)^\frac{d+j}{d}\sqrt{{(1-\gamma)}-
\frac{1}{n\delta_n}}}\Bigg(\sqrt{2}+\sqrt{2}\frac{\kappa_{n}(\alpha)+\eta}
{\Gamma(\delta_n)}-\frac{2\sqrt{2}\log\big((1-\gamma)-\frac{1}{n\delta_n}\big)}{\Gamma(\delta_n)^2}\Bigg).
\end{align}

\begin{enumerate}
\item  If   $H_+^j(f,{K_n}) \geq D\frac{\Gamma(\delta_n)}{\sqrt{2n}}$
 for some constant   $ D\geq\mathbf{D}^j(n,\delta_n)\vee 2(d+j)$, then
\[{\bb{P}}_{\mathcal{N}_n}\left( T_{K_n}\leq c_{K_n}(\alpha)\big|N_n\right)\leq\exp\left(-\frac{\eta^2}{6}\right).\]
\item If   $H_-^j(f,{K_n}) \leq -D\frac{\Gamma(\delta_n)}{\sqrt{2n}}$
 for some constant $ D\geq\mathbf{D}^j(n,\delta_n)\vee  \frac{2d(d+j)}{j}$, then
\[{\bb{P}}_{\mathcal{N}_n}\left( T_{K_n}\geq- c_{K_n}(\alpha)\big|N_n\right)\leq\exp\left(-\frac{\eta^2}{6}\right).\]
\end{enumerate}
\end{lemma}

\begin{proof}[\textbf{Proof}]
We only prove   the first part and define
\[c_{n,+}^j:=\frac{\td{\varphi}_n^{d-1}(P_{e^n}X_{(N_n)})^{d+j}\inf_{K_n} f^{(j)}}{\delta_{N_n}}.\]
Then a tedious but straightforward calculation  shows that the   inequality
\begin{align*}
\frac{2^{d-1}jc_{n,+}^j(N_n-1)}{(2d+j)(d+j)}-\eta\sqrt{\frac{N_n-1}{3}}\geq c_{K_n}(\alpha)
\end{align*}
holds for ${\frac{|K_{N_n}|}{|{K_n}|}\geq1-\ve}$ and $N_n\geq(1-\gamma)n\delta_n$.
This implies
\begin{align*}
{\bb{P}}_{\mathcal{N}_n}\big ( T_{K_n}\leq c_{K_n}(\alpha)\big|N_n\big )
=&~\bb{E} \big [{\bb{P}}_{\mathcal{N}_n}\big ( T_{K_n}\leq c_{K_n}(\alpha)|N_n,X_{(N_n)}\big )\big ] \\
\leq&~\bb{E}\Big [{\bb{P}}_{\mathcal{N}_n}\Big ( T_{K_n}\leq\frac{2^{d-1}jc^j_{n,+}(N_n-1)}{(2d+j)(d+j)}-\eta\sqrt{\frac{N_n-1}{3}}\Big |~N_n,X_{(N_n)}\Big)
\Big]  \\
\leq & ~\bb{E}\Big [ {\bb{P}}_{\mathcal{N}_n}\Big( R_{K_n}
\leq\frac{2^{d-1}jc^j_{n,+}(N_n-1)}{(2d+j)(d+j)}-\eta\sqrt{\frac{N_n-1}{3}} \Big |~N_n,X_{(N_n)}\Big )\Big ] ,
\end{align*}
where we used Lemma \ref{1109}  
   and the notation $R_{K_n} = \sum_{i=1}^{N_n-1}\beta ((\ol{F}_+^j)^{-1}(\td{F}_n(P_{e^n}X_{(i)})))$.
Therefore, the assertion follows from  Lemma \ref{03}  and  Hoeffding's inequality.
\end{proof}

\begin{proof}[\textbf{Proof of Theorem \ref{98}}]
For a proof of the first part we proceed in two steps: firstly, we will find an upper bound for the
probability that the test will not reject for one single wedge. Secondly,   we will consider the probability
for simultaneous rejection on every wedge in $\ca{K}_n$.  For a fixed wedge $K_n\in\ca{K}_n$,
$0<\ve<1$ and $\gamma\in(0,\frac{1}{2}]$, we have 
$$
\bb{P}\big(T_{K_n}\leq c_{K_n}(\alpha)\tn{ for any single }{K_n}\in \ca{K}_n\big)
 \leq{\bb{P}}_{\mathcal{N}_n} \big(T_{K_n}\leq c_{K_n}(\alpha)\tn{ for any single }{K_n}\in \ca{K}_n\big) + \mathbb{P}(\mathcal{N}^c_n), $$
 where the event $\mathcal{N}_n$ is defined in \eqref{Nn}.
Notice that the assumptions of Theorem \ref{98} imply those of Lemma \ref{05} and recall that $F(K_n)=\delta_n$. We have
from the assumption in Theorem \ref{98} \textit{(i)}
 with $K_n=K^i_n$
\[H_+^j(f,{K_n})\geq D\big(1+ {b_n} \tfrac{\sqrt{2}}{\Gamma(\delta_n)}
\big)
\tfrac{\Gamma(\delta_n)}{\sqrt{2n}} ,
\]
which relaxes the assumption on the constant $D$ in Lemma \ref{05} as follows.
Let $\eta_{n,{K_n}}>0$ and
 \begin{align}\label{20}D\big(1+{b_n} \tfrac{\sqrt{2}}{\Gamma(\delta_n)} \big)
 \geq\max\Big\{ \mathbf{D}^j(\eta_{n,{K_n}},\delta_n)
 ,{2(d+j)}{}\Big\},\end{align}
where $\mathbf{D}^j$ is defined in \eqref{D}.
Therefore, it follows from
 Lemma \ref{05}, Lemma \ref{04} and Lemma \ref{07} that the probability under consideration can be bounded   by
\begin{align} \label{-93}
\exp\big (-\tfrac{\eta_{n,{K_n}}^2}{6}\big)+ \exp\Big (-\tfrac{D}{d+j}\sqrt{\tfrac{n\delta_n}{2}}\Gamma(\delta_n) (1-(1-\ve)^\frac{d+j}{d} )\big )
+\exp \big(-n \delta_n\tfrac{\gamma^2}{2}\big ),
\end{align}
which concludes the proof for any  single wedge.

We now consider the union of all wedges of $\ca{K}_n$ and define $\td{\delta}_n:=\inf_{{K_n}\in\ca{K}_n}F(K_n)$. As
the wedges in $\ca{K}_n$ are pairwise disjoint, it follows from Lemma \ref{06}    that
\beq{-95}M_n=\big(\tn{\#}\{{K_n}:\;{K_n}\in\ca{K}_n\}\big)\leq\frac{1}{\td{\delta}_n}=o(n).\eeq
Therefore,  $\ca{K}_n$ consists of a finite number of wedges. From \eqref{-93} and the monotonicity of the
function $\delta\mapsto{\delta\log\big(\frac{\exp(1)}{\delta}\big)}$, we obtain the estimate
\begin{align}
\label{-94}&\bb{P}\big(T_{K_n}\leq c_{K_n}(\alpha)\tn{ for at least one  }{K_n}\in \ca{K}_n\big) \\
& \leq\sum_{{K_n}\in\ca{K}_n}\exp\big (-\tfrac{\eta_{n,{K_n}}^2}{6}\big)
+M_n\Big( \exp\big (-\tfrac{D}{d+j}\sqrt{\tfrac{n\td{\delta}_n}{2}}\Gamma(\td{\delta}_n)
(1-(1-\ve)^\frac{d+j}{d})\big )
+\exp\big(-\tfrac{n\td{\delta}_n \gamma^2}{2}\big)\Big), \nonumber
\end{align}
if condition \eqref{20} is fulfilled for every $K_n\in\ca{K}_n$.
We now show that the right-hand side of \eqref{-94} vanishes as $n\rightarrow\infty$ by investigating the asymptotic
behaviour of every summand. For the first summand, let $\eta_{n,{K_n}}:=\big(6\log\big(\frac{1}{\delta_n}\big)+b_n\big)^\frac{1}{2}$, then
\[\sum_{{K_n}\in\ca{K}_n}\exp\big(-\tfrac{\eta_{n,{K_n}}^2}{6}\big)=\sum_{{K_n}\in\ca{K}_n}\delta_n
\exp\big(-\tfrac{b_n}{6}\big)=o(1),\]
because  $\sum_{{K_n}\in\ca{K}_n}\delta_n \leq 1$ and $b_n\to\infty$ as  $n\to\infty$.
Next, we consider the second summand in \eqref{-94}. An application of Lemma \ref{06} \textit{(b)} gives
\bal\sqrt{\tfrac{n\td{\delta}_n}{2}}\Gamma(\td{\delta}_n)
\geq\sqrt{\td{c}^j_+}\log(\exp(1)n)(1-o(1)).\end{align*}
Hence, by \eqref{-95}, if
\beq{100}\tfrac{D}{d+j}\sqrt{\tfrac{D^2}{(d+j)^2}}\big(1-(1-\ve)^\frac{d+j}{d}\big)>1,\eeq
\bal& M_n \exp\big(-\tfrac{D}{d+j}\sqrt{\td{c}^j_+}\log(\exp(1)n)\big(1-(1-\ve)^\frac{d+j}{d}\big)(1-o(1))\big)
= o(1).\end{align*}

Finally, by \eqref{-95} and Lemma \ref{06} \textit{(b)}, we have
\bal o(1)\;\exp\big(-n\td{\delta}_n\frac{\gamma^2}{2}+\log(n)\big)
&\leq o(1)\;\exp\big(-(1-o(1))\td{c}^j_+\log(\exp(1)n)\frac{\gamma^2}{2}+\log(n)\big)\\
&\leq o(1)\;\exp\big(-\log(n)\big(\td{c}^j_+\frac{\gamma^2}{2}-\big(1+o(1)\big)\big)\big) = o(1),
\end{align*}
if
\beq{o123}\frac{D^2}{(d+j)^2}\frac{\gamma^2}{2}> 1.\eeq
In this case,  the third term vanishes as well as $n\to\infty.$


It remains to show that condition \eqref{20} is fulfilled for every $K_n\in\ca{K}_n$.
With $\kappa_{n}(\alpha)\leq A$, we have to prove that
\begin{align}
\label{21}&D\Big(1+ {b_n}\tfrac{\sqrt{2}}{\Gamma(\delta_n)}\Big)\\
\geq& \frac{ {(2d+j)(d+j)}}{   {\sqrt{3}j} (1-\ve)^\frac{d+j}{d}\sqrt{{(1-\gamma)}-\frac{1}{n\delta_n}}}
\Big(\sqrt{2}+\sqrt{2}\frac{A +\eta_{n,K_n}}{\Gamma(\delta_n) }-\frac{2\sqrt{2}\log\big( 1-\gamma
-\frac{1}{n\delta_n}\big)}{
\Gamma(\delta_n)^2}
 \Big). \nonumber
 \end{align}  From Lemma \ref{06} \textit{(b)} it follows that $n\delta_n\rightarrow\infty$ for $n\rightarrow\infty$ for all
$K_n\in\ca{K}_n$. Thus, using $\eta_{n,{K_n}}:=\big(6\log\big(\frac{1}{\delta_n}\big)+b_n\big)^\frac{1}{2}$,
we find that  for sufficiently large $n$ an upper bound for the right hand side of  \eqref{21} is given by
\begin{align}\begin{split}
\label{-96} & \frac{ {(2d+j)(d+j)}  }{{\sqrt{3}j} (1-\ve)^\frac{d+j}{d}\sqrt{{(1-\gamma)}-o(1)}
}\Big(\sqrt{2}+\frac{A +\big ( \{6\log\big(\tfrac{1}{\delta_n}\big)\}^{1/2}  +\sqrt{b_n}\big )}
{ \big\{\log\big(\tfrac{1}{\delta_n}\big)\big\}^{1/2} }\Big)\\
&\leq\frac{  {(\sqrt{2}+\sqrt{6})(2d+j)(d+j)} }{ { \sqrt{3}j }  (1-\ve)^\frac{d+j}{d}\sqrt{1-\gamma-o(1)}}
+\frac{  {(2d+j)(d+j)}(A+\sqrt{b_n})}{ { \sqrt{3}j }  (1-\ve)^\frac{d+j}{d}
 \big\{\log\big(\tfrac{1}{\delta_n}\big)\big\}^{1/2}
\sqrt{1-\gamma-o(1)}}\\
&\leq\Big(1+\frac{o(b_n)}{  \{\log\big(\tfrac{1}{\delta_n}\big)\}^{1/2}   }\Big)\frac{{2\sqrt{2}(2d+j)(d+j)}{}}
{j(1-\ve)^\frac{d+j}{d}\sqrt{1-\gamma-o(1)}}.
\end{split}
\end{align}
Combining \eqref{100}, \eqref{o123}   and \eqref{-96}, we obtain the following condition
\begin{align}
\frac{D}{d+j}>\max\Big \{\big (1-(1-\ve)^\frac{d+j}{d}\big)^{-{1}/{2}}\,,\,
\tfrac{\sqrt{2}}{\gamma}\,,\,
\tfrac{2\sqrt{2}(2d+j)}{j(1-\ve)^\frac{d+j}{d}\sqrt{1-\gamma}}
\Big \}\label{202}.
\end{align}
In order to minimize the restrictions imposed by condition \eqref{202}, we now determine $0<\ve<1$ and $\gamma\in(0,\frac{1}{2}]$, 
such that the lower bound on $D$  is as small as possible. Balancing the second and third
 terms in \eqref{202} we obtain
 \[
\gamma=\frac{-j^2+\sqrt{j^4+4j^2(2d+j)^2}}{2(2d+j)^2} <\frac{j}{2d+j}\leq\frac{1}{2},
\]
where we used  $(2d+j)^2\geq9$ (note that  $d\geq 1$) for the first inequality.
For the choice of $\ve$ we introduce the notation  $a:=(1-\ve)^\frac{d+j}{d}$
and balance the first and third  expression in \eqref{202} and obtain
\[
a=\tfrac{(2d+j)^2}{j^2}\Big (-1+
\sqrt{1+\tfrac{2j^2}{(2d+j)^2}}\Big ).
\]
Finally, inserting our choice of $\ve$ and $\gamma$ in \eqref{202}, we find the condition \eqref{63i}
since our calculations also show that  $D$ is larger than all three terms of \eqref{202} simultaneously in
this case.
\end{proof}

\subsection{Proof of Theorem \ref{S1A} and  \ref{S1B}}

For the sake of simplicity, we prove both results for the  case  $C_2=1$ and $C:=C_1$.
The general case follows by exactly the same arguments with an additional amount of notation.

\medskip

{\bf Proof of Theorem \ref{S1A}:}  We note that it follows from Theorem \ref{C1} that $c_{K^i}(\alpha)\geq\td{c}_{K^i}(\alpha)$ for $i=1,\hdots,M_n$.
Hence, it remains to show that the assumptions for Theorem \ref{98} \textit{(ii)} are satisfied.
By assumption on $f$, we have
$\sup_{K^i}f^{(2)}\leq -\frac{c}{2}+o(1)\tn{ for }n\rightarrow\infty$ $(i=1,\hdots,M_n)$.
Moreover, from the approximation $\tan(x)=x(1-o(1))$ for $x\rightarrow0$, we have
\beq{-99}F({K}^i)=f(x_0)\frac{1}{d}C^d\log(n)^{-d+1+d\frac{d-1}{d+4}}\Big(\frac{\log(n)}{n}\Big)^{\frac{d}{d+4}}(1+o(1))~~(i=1,\hdots,M_n).\eeq
Hence,
\bal -H_-^2(f,{K^i})&\geq\frac{\log(n)^{-d+1+(d+2)\frac{d-1}{d+4}}C^{d+2}\big(\frac{\log(n)}{n}\big)^
{\frac{d+2}{d+4}}(c/2-o(1))}{\sqrt{f(x_0)\frac{1}{d}}
C^\frac{d}{2}\log(n)^{\frac{-d+1}{2}+d\frac{d-1}{2(d+4)}}\big(\frac{\log(n)}{n}\big)^\frac{d}{2(d+4)}(1+o(1))}
= C^{\frac{d}{2}+2}\big (\tfrac{\log(n)}{n}\big )^{\frac{1}{2}} \frac{c/2-o(1)}{\sqrt{f(x_0)\tfrac{1}{d}}}\end{align*}
$(i=1,\hdots,M_n)$.
Furthermore, we obtain from \eqref{-99} and the assumption $b_n=o(\sqrt{\log(n)})$ that
\bal D\left(\tfrac{\Gamma(F(K^i))}{\sqrt{2n}}+\tfrac{b_n}{\sqrt{n}}\right)
=D\sqrt{\tfrac{d}{d+4}\tfrac{\log n}{n}}(1+o(1)).\end{align*}
Therefore, the assumptions of Theorem \ref{98}  \textit{(ii)} (for $j=2$) are fulfilled as the constant $C$ satisfies $C>\Big(\frac{2D}{c}
\sqrt{\frac{f(x_0)}{d+4}}\Big)^{\frac{2}{d+4}}$ by \eqref{c1}
and hence the assertion  follows by an application of Theorem \ref{98}. \hfill $\Box$ 

\medskip

{\bf Proof of Theorem \ref{S1B}:}
Note that it is sufficient to prove consistency for the largest scale.
By Theorem \ref{C1}, we have
to prove that the assumptions for Theorem \ref{98} \textit{(ii)} are
satisfied for the family of wedges  $\{K^i_n\,|\;i\in \ca{I}_n\}$ introduced in Section \ref{sec33}. Let $K_n\in\{K_n^i\,|\;i\in\ca{I}_n\}$. We begin with the determination of an upper bound for the quantity  $\sup_{K_n} f^{(1)}$  defined in \eqref{supf}.
For this purpose, consider a point  $x_0^n+\td{t}e_0\in K_n$  with  $e_0\in\bb{R}^d$ ($\|e_0\|=1$)
 and $\td{t} >0 $. Now, the representation \eqref{frep}
and an application of  the mean value theorem   yields for $0\leq \td{s}<\td{t}$
\begin{align}
\label{-97}
\Delta_n &= ~\frac{f(x^n_0+\td{t}e_0)-f(x^n_0+\td{s}e_0)}{\td{t}-\td{s}} \\
& =  \td{f}_{x_0}(\|x_0^n+\td{s}e_0-x_0\|)\langle\nabla g_{x_0}(\xi_2),e_0\rangle  
 +(1+g_{x_0}(x^n_0+\td{t}e_0))\td{f}_{x_0}'(\xi_1) R_n  \nonumber
\end{align}
for some $\|x_0^n+\td{s}e_0-x_0\|\leq \xi_1\leq \|x_0^n+\td{t}e_0-x_0\|$ and $ \xi_2\in[x^n_0+\td{s}e_0,x^n_0+\td{t}e_0]$,
where
$$
R_n =  \frac{\|x_0^n+\td{t}e_0-x_0\|
-\|x_0^n+\td{s}e_0-x_0\|}{\td{t}-\td{s}}.
$$
A further application of the mean value theorem gives
 \begin{align*}
 R_n
 &~=\frac{\sum_{j=1}^d(x^n_{0,j}+\td{t}e_{0,j}-x_{0,j})^2-\sum_{j=1}^d(x^n_{0,j}+\td{s}e_{0,j}-x_{0,j})^2}{2\sqrt{\xi}(\td{t}-\td{s})}\\
 &~=\frac{\sum_{j=1}^d2(x^n_{0,j}-x_{0,j})e_{0,j}(\td{t}-\td{s})+\td{t}^2-\td{s}^2}{2\sqrt{\xi}(\td{t}-\td{s})}
 \geq \frac{\sum_{j=1}^d(x^n_{0,j}-x_{0,j})e_{0,j}}{\sqrt{\xi}}\\
 &~=\frac{\cos(\text{angle}(x^n_0-x_0,e_0))\|x^n_0-x_0\|}{\sqrt{\xi}}
 \geq C\log(n)l_n(1-o(1)) \frac{1}{\sqrt{\xi}}
 \end{align*}
  for some $\|x_0^n+\td{s}e_0-x_0\|^2\leq\xi\leq\|x_0^n+\td{t}e_0-x_0\|^2$.
Moreover, we have
\begin{align*}{\xi}\leq\|x^n_0+\td{t}e_0-x_0\|^2\leq\big(\|x^n_0-x_0\|+\td{t}\big)^2
&\leq 
l^2_n(m_nC\log(n))^2(1+o(1)).
 \end{align*}
Hence,
 $ R_n \geq \frac{1-o(1)}{m_n}.$ With the same arguments as used before, one shows that $\xi_1\geq C\log(n)l_n$ and
 $\|\xi_2-x_0\|\leq Cm_n\log(n)l_n(1+o(1))$.
Finally, by assumption on $g_{x_0}$ and $\tilde f_{x_0}$ and \eqref{-97}, this yields
 \begin{align*}
\Delta_n  & ~\leq(1+o(1))(-c)C\log(n)l_n R_n
 +f(x_0)o((m_n\log(n)l_n)^{1+\gamma})\\
&~\leq(1-o(1))(-c)C\log(n)\frac{l_n}{m_n}+o((m_n\log(n)l_n)^{1+\gamma})
=-\frac{cC}{m_n}\log(n)l_n(1-o(1)),\end{align*}
as $ \frac{m_n}{l_n}o((m_nl_n)^{1+\gamma})=o(m_n^{2+\gamma}l_n^\gamma) =o (1)$ as  $n\rightarrow\infty$
by the choice of $m_n$ and $l_n$. Consequently,
$\sup_{K_n}f'\leq-\frac{cC}{m_n}\log(n)l_n(1-o(1)).$
As $\tan(x)=x(1-o(1))$ ($x\rightarrow 0$), we have
\bal F(K_n)&= f(x^n_0)\frac{1}{d}C^d(\log(n))^{-d+1+d}\Big (\frac{\log(n)}{n}\Big )^{\frac{d}{d+4}}(1+o(1))\\
&= f(x_0)\frac{1}{d}C^d\log(n)\Big (\frac{\log(n)}{n}\Big )^{\frac{d}{d+4}}(1+o(1)).\end{align*}
It follows from the conditions  $c_1 \geq  f(x_0)\geq 0$
\bal {H}^1_-(f,{K_n})&\leq-\frac{C^{d+1}(\log(n))^{3}\big (\frac{\log(n)}{n}\big )^{\frac{d+2}{d+4}}
\frac{cC}{m_n}(1-o(1))}{\sqrt{f(x_0)\frac{1}{d}}C^\frac{d}{2}(\log(n))^\frac{1}{2}
\big (\frac{\log(n)}{n}\big )^{\frac{d}{2(d+4)}}(1+o(1))}\\
&\leq -{cC^{\frac{d}{2}+2}}{}\Big (d \frac{\log(n)}{n}\Big )^{\frac{1}{2}}
\frac{1-o(1)}{\sqrt{c_1 }}\frac{(\log(n))^{\frac{5}{2}}}{m_n}
=-{cC^{\frac{d}{2}+2}}{}\Big (d \frac{\log(n)}{n}\Big )^{\frac{1}{2}} \frac{1-o(1)}{\sqrt{c_1  }}.\end{align*}
If $b_n=o(\sqrt{\log n})$, we have
\bal D\Big (\tfrac{\Gamma(F(K_n))}{\sqrt{2n}}+\tfrac{b_n}{\sqrt{n}}\Big)
=D\sqrt{\tfrac{d}{d+4}\tfrac{\log n}{n}}(1+o(1)),
\end{align*}
and the assumptions of Theorem \ref{98} \textit{(ii)} are fulfilled  if the constant $C$ satisfies 
$C > \left( \frac {D}{c}\sqrt{\frac{c_1}{d+4}} \right)^{\frac{2}{d+4}}$, which is a direct consequence of \eqref{bed1}.
\hfill $\Box$

\end{document}

%% file: keil.tex
\psset{xunit=.5pt,yunit=.5pt,runit=.5pt}
\psset{xunit=.5pt,yunit=.5pt,runit=.5pt}
\begin{pspicture}(-350,552.36218262)(-50,280)
{
\newrgbcolor{curcolor}{0 0 0}
\pscustom[linewidth=1,linecolor=curcolor]
{
\newpath
\moveto(-50.5,466.00000262)
\lineto(302.5,594.50000262)
\lineto(302.5,338.50000262)
\closepath
}
}
\psline[linewidth=1,linecolor=curcolor, arrows=->,arrowsize=3pt 4](-50.5,466.00000262)(335.5,466.00000262)
{
\newrgbcolor{curcolor}{0 0 0}
\pscustom[linewidth=1,linecolor=curcolor]
{
\newpath
\moveto(230.5,364.50000262)
\lineto(230.5,568.00000262)
}
}
{
\newrgbcolor{curcolor}{0 0 0}
\pscustom[linewidth=1,linecolor=curcolor]
{
}
}
 \pscurve[linecolor=curcolor,linewidth=1]{}(113.5,466.50000262)(122,481)(125,495.25)(122,510)(113.5,525.50000262)

{
\newrgbcolor{curcolor}{0 0 0}
\pscustom[linewidth=1,linecolor=curcolor]
{
}
}\rput(126,291){$K$}
\rput(170,425){$K_N$}
\rput(343,466.5){$e$}
\rput(310,477.5){$l$}
\rput(260.5,516){$X_{(N)}$}
\rput(230.5,520.00000262){$\times$}
\rput(266.5,448){$P_eX_{(N)}$}
\rput(230.5,466.00000262){$\times$}
\rput(-66,466){$x_0$}
\rput(75, 489){$\vp$}
\psbrace[rot=90,nodesepB=-6pt](-50.5,338.50000262)(302.5,338.50000262){}
\end{pspicture}

%% file: Gitter_nocolor.tex
\psset{xunit=1.2cm,yunit=1.2cm,runit=1.cm}
\begin{pspicture}(-3.5,-2)(3.5,4.5)
\psgrid[subgriddiv=1,%
griddots=10,%
gridlabels=7pt](-4,-2)(4,4)


\psline[linewidth=1.5pt,fillstyle=dots]%
(-3.5,-0.5)(-3,-1)(-2.5,-0.5)(-3.5,-0.5)
\psline[linewidth=1.5pt,fillstyle=dots]%
(-2.5,-0.5)(-3,-1)(-2.5,-1.5)(-2.5,-0.5)

\psline[linewidth=1.5pt,fillstyle=dots]%
(-2.5,-0.5)(-2,-1)(-1.5,-0.5)(-2.5,-0.5)

\psline[linewidth=1.5pt,fillstyle=crosshatch]%
(-1.5,-1.5)(-2,-1)(-2.5,-1.5)(-1.5,-1.5)

\psline[linewidth=1.5pt,fillstyle=dots]%
(-1.5,-1.5)(-1,-1)(-1.5,-0.5)(-1.5,-1.5)
\psline[linewidth=1.5pt,fillstyle=dots]%
(-1.5,-0.5)(-1,-1)(-0.5,-0.5)(-1.5,-0.5)

\psline[linewidth=1.5pt,fillstyle=dots]%
(0.5,-0.5)(1,-1)(1.5,-0.5)(0.5,-0.5)
\psline[linewidth=1.5pt,fillstyle=dots]%
(1.5,-0.5)(1,-1)(1.5,-1.5)(1.5,-0.5)

\psline[linewidth=1.5pt,fillstyle=dots]%
(1.5,-0.5)(2,-1)(2.5,-0.5)(1.5,-0.5)

\psline[linewidth=1.5pt,fillstyle=crosshatch]%
(2.5,-1.5)(2,-1)(1.5,-1.5)(2.5,-1.5)

\psline[linewidth=1.5pt,fillstyle=dots]%
(2.5,-1.5)(3,-1)(2.5,-0.5)(2.5,-1.5)
\psline[linewidth=1.5pt,fillstyle=dots]%
(2.5,-0.5)(3,-1)(3.5,-0.5)(2.5,-0.5)


\psline[linewidth=1.5pt,fillstyle=dots]%
(-2.5,0.5)(-3,0)(-2.5,-0.5)(-2.5,0.5)

\psline[linewidth=1.5pt,fillstyle=crosshatch]%
(-2.5,-0.5)(-2,-0)(-2.5,0.5)(-2.5,-0.5)
\psline[linewidth=1.5pt,fillstyle=crosshatch]%
(-2.5,0.5)(-2,-0)(-1.5,0.5)(-2.5,0.5)
\psline[linewidth=1.5pt,fillstyle=crosshatch]%
(-1.5,0.5)(-2,0)(-1.5,-0.5)(-1.5,0.5)
\psline[linewidth=1.5pt,fillstyle=crosshatch]%
(-1.5,-0.5)(-2,-0)(-2.5,-0.5)(-1.5,-0.5)

\psline[linewidth=1.5pt,fillstyle=dots]%
(-1.5,-0.5)(-1,-0)(-1.5,0.5)(-1.5,0.5)

\psline[linewidth=1.5pt,fillstyle=crosshatch]%
(-0.5,0.5)(-1,0)(-0.5,-0.5)(-0.5,0.5)

\psline[linewidth=1.5pt,fillstyle=crosshatch]%
(0.5,-0.5)(1,-0)(0.5,0.5)(0.5,-0.5)

\psline[linewidth=1.5pt,fillstyle=dots]%
(1.5,0.5)(1,-0)(1.5,-0.5)(1.5,0.5)

\psline[linewidth=1.5pt,fillstyle=crosshatch]%
(1.5,-0.5)(2,-0)(1.5,0.5)(1.5,-0.5)
\psline[linewidth=1.5pt,fillstyle=crosshatch](1.5,0.5)(2,-0)(2.5,0.5)(1.5,0.5)
\psline[linewidth=1.5pt,fillstyle=crosshatch](2.5,0.5)(2,-0)(2.5,-0.5)(2.5,0.5)
\psline[linewidth=1.5pt,fillstyle=crosshatch](2.5,-0.5)(2,-0)(1.5,-0.5)(2.5,-0.5)

\psline[linewidth=1.5pt,fillstyle=dots]%
(2.5,-0.5)(3,-0)(2.5,0.5)(2.5,-0.5)

\psline[linewidth=1.5pt,fillstyle=crosshatch]%
(3.5,0.5)(3,-0)(3.5,-0.5)(3.5,0.5)


\psline[linewidth=1.5pt,fillstyle=dots]%
(-1.5,0.5)(-2,1)(-2.5,0.5)(-1.5,0.5)

\psline[linewidth=1.5pt,fillstyle=dots](2.5,0.5)(2,1)(1.5,0.5)(2.5,0.5)


\psline[linewidth=1.5pt,fillstyle=crosshatch]%
(-1.5,1.5)(-1,2)(-1.5,2.5)(-1.5,1.5)

\psline[linewidth=1.5pt,fillstyle=dots]%
(-0.5,2.5)(-1,2)(-0.5,1.5)(-0.5,2.5)

\psline[linewidth=1.5pt,fillstyle=crosshatch]%
(0.5,2.5)(0,2)(0.5,1.5)(0.5,2.5)
\psline[linewidth=1.5pt,fillstyle=crosshatch](0.5,1.5)(0,2)(-0.5,1.5)(0.5,1.5)

\psline[linewidth=1.5pt,fillstyle=crosshatch]%
(1.5,2.5)(1,2)(1.5,1.5)(1.5,2.5)


\psline[linewidth=1.5pt,fillstyle=dots]%
(0.5,2.5)(0,3)(-0.5,2.5)(0.5,2.5)

\psline[linewidth=0.5pt](-3.5,-1.5)(3.5,-1.5)
\psline[linewidth=0.5pt](-3.5,-0.5)(3.5,-0.5)
\psline[linewidth=0.5pt](-3.5,0.5)(3.5,0.5)
\psline[linewidth=0.5pt](-3.5,1.5)(3.5,1.5)
\psline[linewidth=0.5pt](-3.5,2.5)(3.5,2.5)
\psline[linewidth=0.5pt](-3.5,2.5)(3.5,2.5)
\psline[linewidth=0.5pt](-3.5,3.5)(3.5,3.5)

\psline[linewidth=0.5pt](-3.5,-1.5)(-3.5,3.5)
\psline[linewidth=0.5pt](-2.5,-1.5)(-2.5,3.5)
\psline[linewidth=0.5pt](-1.5,-1.5)(-1.5,3.5)
\psline[linewidth=0.5pt](-0.5,-1.5)(-0.5,3.5)
\psline[linewidth=0.5pt](0.5,-1.5)(0.5,3.5)
\psline[linewidth=0.5pt](1.5,-1.5)(1.5,3.5)
\psline[linewidth=0.5pt](2.5,-1.5)(2.5,3.5)
\psline[linewidth=0.5pt](3.5,-1.5)(3.5,3.5)

\psline[linewidth=0.5pt](-3.5,-1.5)(1.5,3.5)
\psline[linewidth=0.5pt](-2.5,-1.5)(2.5,3.5)
\psline[linewidth=0.5pt](-1.5,-1.5)(3.5,3.5)
\psline[linewidth=0.5pt](-1.5,-1.5)(3.5,3.5)
\psline[linewidth=0.5pt](-0.5,-1.5)(3.5,2.5)
\psline[linewidth=0.5pt](0.5,-1.5)(3.5,1.5)
\psline[linewidth=0.5pt](1.5,-1.5)(3.5,0.5)
\psline[linewidth=0.5pt](2.5,-1.5)(3.5,-0.5)
\psline[linewidth=0.5pt](-3.5,-0.5)(0.5,3.5)
\psline[linewidth=0.5pt](-3.5,0.5)(-0.5,3.5)
\psline[linewidth=0.5pt](-3.5,1.5)(-1.5,3.5)
\psline[linewidth=0.5pt](-3.5,2.5)(-2.5,3.5)

\psline[linewidth=0.5pt](-3.5,-0.5)(-2.5,-1.5)
\psline[linewidth=0.5pt](-3.5,0.5)(-1.5,-1.5)
\psline[linewidth=0.5pt](-3.5,1.5)(-0.5,-1.5)
\psline[linewidth=0.5pt](-3.5,2.5)(0.5,-1.5)
\psline[linewidth=0.5pt](-3.5,3.5)(1.5,-1.5)
\psline[linewidth=0.5pt](-2.5,3.5)(2.5,-1.5)
\psline[linewidth=0.5pt](-1.5,3.5)(3.5,-1.5)
\psline[linewidth=0.5pt](-0.5,3.5)(3.5,-0.5)
\psline[linewidth=0.5pt](0.5,3.5)(3.5,0.5)
\psline[linewidth=0.5pt](1.5,3.5)(3.5,1.5)
\psline[linewidth=0.5pt](2.5,3.5)(3.5,2.5)

\end{pspicture}

%% file: Keil2.tex
\psset{xunit=.5pt,yunit=.5pt,runit=.5pt}
\psset{xunit=.5pt,yunit=.5pt,runit=.5pt}
\begin{pspicture}(-350,552.36218262)(-50,270)
{
\newrgbcolor{curcolor}{0 0 0}
\pscustom[linewidth=1,linecolor=curcolor]
{
\newpath
\moveto(-50.5,466.00000262)
\lineto(302.5,594.50000262)
\lineto(302.5,338.50000262)
\closepath
}
}
\psline[linewidth=1,linecolor=curcolor, arrows=->,arrowsize=3pt 4](-50.5,466.00000262)(335.5,466.00000262)
{
\newrgbcolor{curcolor}{0 0 0}
\pscustom[linewidth=1,linecolor=curcolor]
{
\newpath
\moveto(230.5,364.50000262)
\lineto(230.5,568.00000262)
}
}
{
\newrgbcolor{curcolor}{0 0 0}
\pscustom[linewidth=1,linecolor=curcolor]
{
\newpath
\moveto(105.5,410.0000262)
\lineto(105.5,522.50000262)
}
}
{
\newrgbcolor{curcolor}{0 0 0}
\pscustom[linewidth=1,linecolor=curcolor]
{
}
}
\psbrace[rot=90,nodesepB=-6pt](105.5,364.50000262)(230.5,364.50000262){}

{
\newrgbcolor{curcolor}{0 0 0}
\pscustom[linewidth=1,linecolor=curcolor]
{
}
}\rput(167,317){$K(j,k)$}
\rput(126,265){$K$}
\rput(343,466.5){$e$}
\rput(310,477.5){$l$}
\rput(-66,466){$x_0$}
\rput(256.5,376){$X_{(k)}$}
\rput(230.5,380.00000262){$\times$}
\rput(131.5,496){$X_{(j)}$}
\rput(105.5,500.00000262){$\times$}
\rput(264.5,448){$PX_{(k)}$}
\rput(230.5,466.00000262){$\times$}
\rput(139.5,448){$PX_{(j)}$}
\rput(105.5,466.00000262){$\times$}

\psbrace[rot=90,nodesepB=-6pt](-50.5,313.50000262)(302.5,313.50000262){}
\end{pspicture}

%% file: Mode2_nocolor.tex
\psset{xunit=1.2cm,yunit=1.2cm,runit=1.cm}
\begin{pspicture}(2,2.5)(7,4.)

\psline[linewidth=1.5pt,fillstyle=crosshatch]%
(2.5,2.5)(3,3)(2.5,3.5)(2.5,2.5)
\psline[linewidth=1.5pt,fillstyle=crosshatch]%
(2.5,3.5)(3,3)(3.5,3.5)(2.5,3.5)
\psline[linewidth=1.5pt,fillstyle=crosshatch]%
(3.5,3.5)(3,3)(3.5,2.5)(3.5,3.5)
\psline[linewidth=1.5pt,fillstyle=crosshatch]%
(3.5,2.5)(3,3)(2.5,2.5)(3.5,2.5)

\psline[linewidth=1.5pt,fillstyle=dots]%
(2.5,2.5)(2,3)(2.5,3.5)(2.5,2.5)
\psline[linewidth=1.5pt,fillstyle=dots]%
(2.5,2.5)(3,2)(3.5,2.5)(2.5,2.5)
\psline[linewidth=1.5pt,fillstyle=dots]%
(3.5,2.5)(4,3)(3.5,3.5)(3.5,2.5)
\psline[linewidth=1.5pt,fillstyle=dots]%
(3.5,3.5)(3,4)(2.5,3.5)(3.5,3.5)

\psline[linewidth=0.5pt](2.5,2.5)(3.5,3.5)
\psline[linewidth=0.5pt](3,2)(4,3)
\psline[linewidth=0.5pt](2,3)(3,4)
\psline[linewidth=0.5pt](2.5,2.5)(3.5,2.5)
\psline[linewidth=0.5pt](2.5,3.5)(3.5,3.5)
\psline[linewidth=0.5pt](2.5,2.5)(2.5,3.5)
\psline[linewidth=0.5pt](3.5,2.5)(3.5,3.5)
\psline[linewidth=0.5pt](2,3)(3,2)
\psline[linewidth=0.5pt](2.5,3.5)(3.5,2.5)
\psline[linewidth=0.5pt](3,4)(4,3)

\psline[linewidth=1.5pt,fillstyle=crosshatch]%
(6.5,3.5)(6,3)(6.5,2.5)(6.5,3.5)
\psline[linewidth=1.5pt,fillstyle=crosshatch]%
(6.5,2.5)(6,3)(5.5,2.5)(6.5,2.5)

\psline[linewidth=1.5pt,fillstyle=dots]%
(5.5,2.5)(5,3)(5.5,3.5)(5.5,2.5)
\psline[linewidth=1.5pt,fillstyle=dots]%
(6.5,3.5)(6,4)(5.5,3.5)(6.5,3.5)

\psline[linewidth=0.5pt](5.5,2.5)(6.5,3.5)
\psline[linewidth=0.5pt](5,3)(6,4)
\psline[linewidth=0.5pt](5.5,2.5)(6.5,2.5)
\psline[linewidth=0.5pt](5.5,3.5)(6.5,3.5)
\psline[linewidth=0.5pt](5.5,2.5)(6.5,3.5)
\psline[linewidth=0.5pt](6.5,2.5)(6.5,3.5)
\psline[linewidth=0.5pt](5.5,3.5)(6.5,2.5)

\end{pspicture}

%% file: exa117.tex
\begin{center}
\begin{figure}[ht]
\psset{xunit=.5pt,yunit=.5pt,runit=.5pt}
\psset{xunit=.5pt,yunit=.5pt,runit=.5pt}
\begin{pspicture}(-350,607.36218262)(-50,345)
{
\newrgbcolor{curcolor}{0 0 0}
\pscustom[linewidth=1,linecolor=curcolor]
{
\newpath
\moveto(-50.5,466.00000262)
\lineto(302.5,594.50000262)
\lineto(302.5,338.50000262)
\closepath
}
}
{
\newrgbcolor{curcolor}{0 0 0}
\pscustom[linewidth=1,linecolor=curcolor]
{
}
}
\psline[linewidth=1,linecolor=curcolor, arrows=->,arrowsize=3pt 4](-50.5,466.00000262)(322.5,550.50000262)
\psline[linewidth=1,linecolor=curcolor, arrows=->,arrowsize=3pt 4](-50.5,466.00000262)(335.5,466.00000262)
{
\newrgbcolor{curcolor}{0 0 0}
\pscustom[linewidth=1,linecolor=curcolor]
{
\newpath
\moveto(221.5,466.50000262)
\lineto(221.5,527.00000262)
}
}
{
\newrgbcolor{curcolor}{0 0 0}
\pscustom[linewidth=1,linecolor=curcolor]
{
}
}

{
\newrgbcolor{curcolor}{0 0 0}
\pscustom[linewidth=1,linecolor=curcolor]
{
}
}\rput(248,410){$K$}
\rput(343,466.5){$e$}
\rput(334,550.5){$e_0$}
\rput(221.5,452){${t_1}$}
\rput(221.5,540){$x$}
\rput(221.5,527.00000262){$\times$}
\psbrace[rot=180,nodesepB=-6pt](221.5,466.50000262)(221.5,527.00000262){}
\psbrace[rot=180,nodesepB=-6pt](221.5,549.00000262)(-50.5,488.00000262){}
\rput(78,566){$\td{t}_1$}

\rput(268,494){$s_1$}
\rput(-66,466){$x_0$}

\end{pspicture}
\caption{Representation of $x\in K$ for $d=2$.}\label{AA2}
\end{figure}
\end{center}

%% file: irn.tex

\psset{xunit=.5pt,yunit=.5pt,runit=.5pt}
\psset{xunit=.5pt,yunit=.5pt,runit=.5pt}
\begin{pspicture}(-350,552.36218262)(-50,288)
{
\newrgbcolor{curcolor}{0 0 0}
\pscustom[linewidth=1,linecolor=curcolor]
{
\newpath
\moveto(-50.5,466.00000262)
\lineto(250.5,574.00000262)
\lineto(250.5,358.50000262)
\closepath
}
}
{
\newrgbcolor{curcolor}{0 0 0}
\pscustom[linewidth=1,linecolor=curcolor]
{
}
}
{
\newrgbcolor{curcolor}{0 0 0}
\pscustom[linewidth=1,linecolor=curcolor]
{
\newpath
\moveto(250.5,358.50000262)
\lineto(250.5,574.00000262)
}
}
{
\newrgbcolor{curcolor}{0 0 0}
\pscustom[linewidth=1,linecolor=curcolor]
{
}
}

{
\newrgbcolor{curcolor}{0 0 0}
\pscustom[linewidth=1,linecolor=curcolor]
{
}
}\rput(281,288){$I_{r_n}$}
\rput(200,510){$K_n$}
\rput(347,466.5){$e^n$}
\rput(316,450.5){$l_n$}
\rput(107.5,420){$l_n(1-\ve)^\frac{1}{d}$}
\rput(-66,466){$x^n_0$}
\psbrace[rot=90,nodesepB=-6pt](250,338.50000262)(302.5,338.50000262){}
\psline[linewidth=1,linecolor=curcolor, arrows=->,arrowsize=3pt 4](-50.5,466.00000262)(335.5,466.00000262)
\psbrace[rot=90,nodesepB=-6pt](-50.5,466.00000262)(250,466.00000262){}

   \pscustom[fillstyle=vlines]{%
   \psline(250.5,358.50000262)(302.5,338.50000262)(302.5,594.50000262)(250.5,574.00000262)(250.5,358.50000262)
}
\end{pspicture}
